\documentclass{article}
\usepackage{latexsym,amssymb,amsmath}
\usepackage{graphicx}
\usepackage{color}
\usepackage[emtex,matrix,arrow]{xy}
\newtheorem{thm}{Theorem}[section]
\newtheorem{cor}[thm]{Corol$\cdot$lary}
\newtheorem{ex}[thm]{Example}
\newtheorem{lem}[thm]{Lemma}
\newtheorem{prop}[thm]{Proposition}
\newtheorem{defn}[thm]{Definition}

\numberwithin{equation}{section}
\newcommand{\Real}{\mathbb R}

\newcommand{\Complex}{\mathbb C}

\newcommand{\Natural}{\mathbb N}

\newcommand{\To}{\longrightarrow}

\newcommand{\Aa}{\mathcal{A}}

\newcommand{\Bb}{\mathcal{B}}

\newcommand{\Rr}{\mathcal{R}}
\newcommand{\Ss}{\mathcal{S}}

\newcommand{\Ii}{\mathcal{I}}
\newcommand{\Oo}{\mathcal{O}}
\newcommand{\F}{\mathcal{F}}

\newcommand{\G}{\mathcal{G}}

\newcommand{\dosv}{\bf{2Vect}}
\newcommand{\dossv}{\bf{2SVect}}
\newcommand{\PF}{\bf{2Fun}}
\newcommand{\SVe}{\bf{SVect}}
\newcommand{\Ve}{\bf{Vect}}

\newcommand{\repg}{\mathfrak{Rep}}

\newcommand{\Cgg}{\mathfrak{C}}
\newcommand{\ngg}{\mathfrak{n}}
\newcommand{\Ggg}{\mathfrak{G}}
\newcommand{\Dgg}{\mathfrak{D}}
\newcommand{\id}{\rm{id}}

\newcommand{\qed}{\hspace*{\fill}$\Box$  \ifmmode \else
    \par\addvspace\topsep\fi}
\newenvironment {proof}{\par\addvspace\topsep\noindent{\it Proof.}
    \ignorespaces }{\qed}

\begin{document}

\title{A strict totally coordinatized version of Kapranov and
  Voevodsky's 2-category {\dosv}}
\author{Josep Elgueta \\ Dept. Matem\`atica Aplicada II \\ Universitat
  Polit\`ecnica de Catalunya \\ email: Josep.Elgueta@upc.es}

\date{}


\maketitle

\begin{abstract}
The purpose of this paper is to give a concrete description of a
strict totally coordinatized version of Kapranov and Voevodsky's
2-category of finite dimensional 2-vector spaces. In particular, we
give explicit formulas for the composition of 1-morphisms and
the two compositions between 2-morphisms. 
\end{abstract}

\section{Introduction}

Let {\Ve} be the category of finite dimensional vector
spaces over $\Complex$. In \cite{KV94}, Kapranov and Voevodsky
introduced a sort of categorification of {\Ve}, denoted {\dosv} and
called the 2-category of (finite dimensional) {\sl 2-vector
spaces}. Roughly, the idea is to take the (ring) category {\Ve} as
the analog of the field $\Complex$. Thus, objects in {\dosv} are
what they call {\Ve}-module categories {\Ve}-module equivalent to
some finite product ${\Ve}^n$ of {\Ve} (analog of the
$\Complex$-modules linearly isomorphic to $\Complex^n$) and
1-morphisms are the so-called {\Ve}-module functors (analogous to the
$\Complex$-linear maps); cf. \cite{KV94} for more details. Needless
to say, in {\dosv} we additionally have suitably defined
2-morphisms having no analog in {\Ve}.

It is well-known that {\Ve} is equivalent to the skeletal category
{\bf Mat} whose objects are all integers $n\geq 0$ and where each
$m\times n$ complex matrix $A$ ($m,n\geq 1$) is regarded as a morphism $A:n\To
m$, with composition the usual matrix product. The passage from {\Ve} to {\bf
  Mat} involves two steps: a 
``skeletization'' step, consisting in the replacement of {\Ve}
by a skeleton, and a ``coordinatization'' step, where linear
maps are identified with complex matrices by choosing a basis in
each vector space.

Similarly, Kapranov and Voevodsky introduced two
coordinatized versions of suitable skeleta of {\dosv}, denoted
${\dosv}_c$ and ${\dosv}_{cc}$ by the authors. They differ in both
the level of skeletization and the level of coordinatization. Thus,
while in both versions the objects are natural numbers,
corresponding to the fact that {\dosv} is replaced by a skeleton
with the same set of objects in both cases, the description of the
1- and 2-morphisms is different. More explicitly, 1-morphisms in
${\dosv}_c$ are matrices of (finite dimensional) vector spaces,
while in ${\dosv}_{cc}$ they are matrices with entries in the set
$\Natural$ of natural numbers. This is due to an additional
skeletization (not coordinatization) step carried out in
${\dosv}_{cc}$, in which every vector space is identified with its
dimension. As regards 2-morphisms, they are matrices whose entries
are linear maps between corresponding vector spaces in ${\dosv}_c$,
while they are matrices of usual complex matrices in
${\dosv}_{cc}$. The difference is due now to an additional
coordinatization step carried out in ${\dosv}_{cc}$.

As pointed out by the authors, ${\dosv}_c$ is not a {\it strict}
2-category, because composition of 1-morphisms is not strictly
associative. However, contrary to what authors claim  (cf. \cite{KV94},
Definition 5.9), ${\dosv}_{cc}$ is also 
non strict. Thus, composition of 1-morphisms is given by the usual product of
matrices, which is indeed strictly
associative. But this is no longer true for the horizontal
composition of 2-morphisms, which, althought not explicitly stated, seems to
be defined by the usual formula for multiplying matrices but with the product
and sum of
the entries respectively given by the tensor product and the direct sum of
complex matrices. This gives a non strictly associative
composition between 2-morphisms precisely 
because the direct sum of matrices is
non commutative (actually, it seems that Proposition 5.8. in \cite{KV94}, from
which the authors would deduce the strict character of ${\dosv}_{cc}$, will
not be true in general). This implies that the associator can not be taken
trivial, because otherwise identity 2-morphisms would fail to satisfy the
naturality condition required to any associator.

What happens is that,
when the 1-morphisms of {\dosv} are identified with matrices of vector
spaces (skeletization step on morphisms
done to go from {\dosv} to ${\dosv}_c$), the strict associativity of the
composition of 
1-morphisms is lost and, consequently, also the associativity of the
horizontal composition of 2-morphisms. When one further skeletizes the
category of morphisms to get matrices of positive integers as
1-morphisms (step from ${\dosv}_c$ to ${\dosv}_{cc}$), strict
associativity of the composition of 1-morphisms is 
recovered, but not the associativity of the horizontal composition of
2-morphisms. 

The purpose of this paper is to describe a really {\it strict} and totally
coordinatized version of {\dosv}, which will be denoted by
${\dossv}_{cc}$. Since a coordinatization process
clearly preserves the
strict character of {\dosv}, the non
strict character of ${\dosv}_{cc}$ arises in the skeletization steps
carried 
out by the authors on the categories of morphisms. The idea to get our version
${\dossv}_{cc}$ is then very simple, although not completely
obvious to implement it. Namely, to take a suitable skeleton {\dossv}
of {\dosv} which still is a strict
2-category (this is accomplished by skeletizing only at the level of
objects) and to carry out {\it uniquely} a coordinatization process
in order to keep strictness. The result is a 2-category
where the objects are natural numbers, as in ${\dosv}_{cc}$, but where the
1- and 2-morphisms, as well as the various compositions, have a
slightly more involved description than in Kapranov and Voevodsky's
version.

Having at our disposal a strict and manageable version
of {\dosv} seems to be a desirable goal, if one wishes to undertake
a generalization of classical linear algebra and its applications to
the categorical setting and if {\dosv} is really the good candidate for
playing the role of {\Ve}. In this sense, it is worth pointing out the repeated
use various mathematicians have made of {\dosv} (or some version of
it) in the last 
years. Thus, we have the above
mentioned work by Kapranov and Voevodsky \cite{KV94}, where the 2-category
{\dosv} is introduced in order to get what it looks as
the right conceptual framework for Zamolodchikov tetrahedra equations.
There is also
(unpublished) work by Yetter \cite{dYnp} claiming
for a {\it categorical linear algebra} theory as a good candidate to provide 
us with ``an algebraic footing for 
the extension to higher dimensions of the successful interaction between
3-manifold topology, quantum field theory and monoidal category theory''. Along
these lines, Mackaay \cite{mM00}, starting
with a suitable pair of finite groups, introduces some sort of
generalization to the 2-categorical setting of the classical group
algebras, and shows how to 
use them to define invariants of 4-manifolds. In a
different direction, but also related to these developpments, 
Neuchl \cite{mN97} considers the  
representation theory on 2-vector spaces of {\it Hopf categories}, an
analog in the categorical setting of the classical notion of Hopf
algebra, and which was first introduced by Crane and Frenkel
\cite{CF94} as an algebraic input to construct four dimensional
topological quantum field theories. In particular, in 
Neuchl's work the reader may 
find an entire chapter devoted to ``2-dimensional linear 
algebra'' where a more general notion of (possibly
infinite dimensional) 2-vector space is given. More 
recently, Barrett and Mackaay \cite{BMnp}, in unpublished work, started
exploring the 
representation theory of categorical groups on {\dosv} as the natural analog
in the category setting of classical representation theory of groups.

This has actually been our motivation in undertaking this 
work: the representation theory of categorical groups. We start pursuing this
direction in
\cite{jE4}, where we use the 2-category
${\dossv}_{cc}$ to give a very explicit description of the
2-category structure of the representations of any categorical group as
automorphisms of a finite-dimensional 2-vector space. In particular,
it easily follows from this description that the 
representation theory of categorical groups on these objects includes
classical 
representation theory of groups on finite dimensional vector spaces,
in the sense that the monoidal category of 
representations of any group $G$ can be recovered as a full subcategory of
the category of endomorphisms of a particular representation 
of $G$, when $G$ is thought of as a categorical group with only
identity arrows (as it should be suspected, the category of
representations of a categorical group is not a category but a 2-category).

The outline of the paper is as follows. In a very short Section 2 we introduce
the above mentioned skeleton of {\dosv} prior to the coordinatization
process. The complete coordinatization of this skeleton is then
carried out in Section 3, which 
constitutes the core of the paper. In
Section 4 we present the final definition of ${\dossv}_{cc}$.

{\bf Notation}. For any $n\geq 1$, we denote by ${\bf
e}_1,\ldots,{\bf e}_n\in\Natural^n$ the canonical basis of
$\Complex^n$. The matrix in canonical basis of a linear map
$f:\Complex^a\To\Complex^b$ is denoted by $[f]$. For any
$m\times n$ complex matrix {\bf M}, ${\bf M}_i$
$i=1,\ldots,m$, is the vector in $\Complex^m$ defined by the $i^{th}$
row of {\bf M}. Furthermore, if ${\bf x}\in\Complex^n$ $, {\bf
M}({\bf x})$ stands for the usual action of {\bf M} on {\bf x}. In
particular, ${\bf M}({\bf e}_j)$ is the vector in $\Complex^m$
defined by the $j^{th}$ column of {\bf M}. Given points
${\bf a},{\bf r}\in\Natural^n$, ${\bf a\cdot r}$ stands for the
usual dot product.

We denote by ${\SVe}$ the full subcategory of {\Ve} with objects the
spaces $\Complex^n$, for all $n\geq 0$ (the ``S'' stands for
``skeleton''). {\SVe}, and hence also its finite products
${\SVe}^n$ for all $n\geq 1$, are $\Complex$-linear additive
categories. In the whole paper, we will take as biproduct functor
$\oplus=\oplus_{{\SVe}}$ on ${\SVe}$ that acting on morphisms in
the standard way (on objects
it is necessarily given by
$\Complex^a\oplus\Complex^{b}=\Complex^{a+b}$). Thus, given linear
maps $f:\Complex^a\To\Complex^{a'}$ and
$g:\Complex^b\To\Complex^{b'}$, with matrices in canonical basis
$[f]$, $[g]$ respectively, $f\oplus
g:\Complex^{a+b}\To\Complex^{a'+b'}$ denotes the linear map such that
$$
[f\oplus g]=\left(\begin{array}{cc} [f] & 0 \\ 0 & [g]
\end{array}\right)
$$
This corresponds to taking as morphisms
$\iota^a_{a,b}:\Complex^a\To\Complex^{a+b}$,
$\iota^b_{a,b}:\Complex^b\To\Complex^{a+b}$,
$\pi^a_{a,b}:\Complex^{a+b}\To\Complex^a$ and
$\pi^b_{a,b}:\Complex^{a+b}\To\Complex^b$ defining the biproduct of
$\Complex^a$ and $\Complex^b$ the linear maps defined by
\begin{align*}
[\iota^a_{a,b}]&=[\pi^a_{a,b}]^T=\left(\begin{array}{cc} Id_a
\\ {\bf 0}_{b\times a}\end{array} \right) \\
[\iota^b_{a,b}]&=[\pi^b_{a,b}]^T=\left(\begin{array}{cc} {\bf
0}_{b\times a} \\ Id_b \end{array} \right)
\end{align*}
where ${\bf 0}_{p\times q}$ denotes the $p\times q$ zero matrix. As
biproduct functor in ${\SVe}^n$, denoted $\oplus_n$, we 
take that induced by $\oplus$ componentwise. In particular, for any
two objects $(\Complex^{a_1},\ldots,\Complex^{a_n})$ and
$(\Complex^{b_1},\ldots,\Complex^{b_n})$ of ${\SVe}^n$, the canonical maps
\begin{align*}
\iota^{(x_1,\ldots,x_n)}_{(a_1,\ldots,a_n),(b_1,\ldots,b_n)}&:
(\Complex^{x_1},\ldots,\Complex^{x_n})\To
(\Complex^{a_1},\ldots,\Complex^{a_n})\oplus_n
(\Complex^{b_1},\ldots,\Complex^{b_n}) \\ 
\pi^{(x_1,\ldots,x_n)}_{(a_1,\ldots,a_n),(b_1,\ldots,b_n)}&:
(\Complex^{a_1},\ldots,\Complex^{a_n})\oplus_n
(\Complex^{b_1},\ldots,\Complex^{b_n})\To
(\Complex^{x_1},\ldots,\Complex^{x_n})
\end{align*}
are given by
\begin{align*}
\iota^{(x_1,\ldots,x_n)}_{(a_1,\ldots,a_n),(b_1,\ldots,b_n)}
&=(\iota^{x_1}_{a_1,b_1},\ldots,\iota^{x_n}_{a_n,b_n}) \\
\pi^{(x_1,\ldots,x_n)}_{(a_1,\ldots,a_n),(b_1,\ldots,b_n)}
&=(\pi^{x_1}_{a_1,b_1},\ldots,\pi^{x_n}_{a_n,b_n})
\end{align*}
for $x=a,b$, with 
$\iota^{x_j}_{a_j,b_j}:\Complex^{x_j}\To\Complex^{a_j+b_j}$ and
$\pi^{x_j}_{a_j,b_j}:\Complex^{a_j+b_j}\To\Complex^{x_j}$ the maps
given above.

\section{The strict 2-category ${\dossv}$}

To shorten the path from {\dosv} to ${\dossv}_{cc}$, we will take
as our starting point the following 2-category {\dossv}, which is
skeletal on objects, but not on morphisms.

\begin{defn}
Let {\dossv} be the 2-category whose objects, 1-morphisms and
2-morphisms are respectively the $\Complex$-linear additive
categories ${\SVe}^n$ for all $n\geq 0$, the $\Complex$-linear
(hence, additive) functors $F:{\SVe}^n\To{\SVe}^m$ and the natural
transformations, with the usual compositions of functors and
natural transformations and usual 1- and 2-identities.
\end{defn}
In this definition, ${\SVe}^0$ means the terminal category with
only one object and one (identity) morphism. It has an obvious
$\Complex$-linear additive structure.

As it will become clear
in the sequel, choosing as objects of {\dosv} the categories
${\SVe}^n$ instead of ${\Ve}^n$ or any category $\Complex$-linear
equivalent to ${\Ve}^n$, for $n\geq 0$, is an essential step in
getting our strict totally coordinatized version of {\dosv}. It allows
us to identify the objects of ${\SVe}^n$ with points 
${\bf a}=(a_1,\ldots,a_n)\in\Natural^n$ and, most importantly, it
makes possible to completely coordinatize the 2-category in a
relatively easy way without carrying out additional skeletizations.

The reader may be wondering if the 1- and 2-morphisms considered in
the previous definition really correspond to those considered by
Kapranov and Voevodsky, namely, the {\Ve}-module functors and the {\Ve}-module
natural transformations, respectively. This was indeed shown by
Yetter \cite{dYnp}, who proved that the
notion of {\Ve}-module functor between suitable {\Ve}-module categories,
of which our categories ${\SVe}^n$ are particular cases, is equivalent
to that of a $\Complex$-linear functor and, furthermore, that every natural
transformation between such functors is a {\Ve}-module natural
transformation.

\section{Coordinatizing {\dossv}}

As indicated in the introduction, to get our 2-category
${\dossv}_{cc}$ from {\dossv} we will just carry out a
coordinatization process, implying no additional ``skeletization''.
In other words, we want to identify manageable sets parametrizing both
the sets of 1-morphisms between two given objects ${\SVe}^n$ and
${\SVe}^m$ and the sets of 2-morphisms between any two 1-morphisms. As
it will be seen in the sequel, this requires fixing various 
functors (in coordinatizing both the notion of 1- and 2-morphism). This is to
be thought of as an analog of the
choice of a basis in each vector space when going from {\SVe} to
{\bf Mat}.

\subsection{1-morphisms}

Given $n,m\geq 1$, let us denote by ${\bf
Fun}_{\Complex}({\SVe}^n,{\SVe}^m)$ the set of $\Complex$-linear
functors $F:{\SVe}^n\To{\SVe}^m$ and let ${\rm Mat}_{m\times
n}(\Natural)$ be the set of $m\times n$ matrices with entries in
$\Natural$. Recall that a $\Complex$-linear functor
$F:{\SVe}^n\To{\SVe}^m$, for any $m\geq 1$, is the same thing as a
collection of $\Complex$-linear functors
$F_1,\ldots,F_m:{\SVe}^n\To{\SVe}$, called its
{\it components}.

For any $j=1,\ldots,n$, let us denote by $\Complex(j,n)$ the object of
${\SVe}^n$ given by
$$
\Complex(j,n)=(0,\ldots,0,\stackrel{j)}{\Complex},0,\ldots,0),\quad j=1,\ldots,n
$$
Then, we have a map ${\bf R}:{\bf
Fun}_{\Complex}({\SVe}^n,{\SVe}^m)\To{\rm Mat}_{m\times
n}(\Natural)$ which applies a $\Complex$-linear functor $F$ to the
matrix ${\bf R}^F$ defined by
\begin{equation} \label{definicio_matriu_rangs}
{\bf R}^F_{ij}={\rm dim}\ F_i(\Complex(j,n)),\quad i=1,\ldots,m,\
j=1,\ldots,n
\end{equation}
where $F_i$ is the $i^{th}$-component of $F$, $i=1,\ldots,m$.
${\bf R}^F$ will be called the {\sl rank matrix} (or {\sl rank
vector} in case $m=1$) of $F$. Thus, we have
\begin{equation} \label{accio_sobre_C(j,n)}
F(\Complex(j,n))=(\Complex^{{\bf R}^F_{1j}},\ldots,\Complex^{{\bf
R}^F_{mj}})\quad j=1,\ldots,n
\end{equation}
This matrix is the analog in our setting of the $m\times n$ complex
matrix describing a morphism $f:\Complex^n\To\Complex^m$ in {\Ve}.
The choice of a basis in the domain and codomain vector spaces
$\Complex^n$ and $\Complex^m$ corresponds in our case to the choice
of the objects $\{\Complex(j,n),\ j=1,\ldots,n\}$ and
$\{\Complex(i,m),\ i=1,\ldots,m\}$ as the basic ones from which any
other object may be constructed by taking finite biproducts.
In the same way as the matrix of the linear map
completely determines its action on the vectors of $\Complex^n$,
${\bf R}^F$ uniquely determines the action of $F$ on objects.
More explicitly, any $\Complex$-linear functor is in particular additive and
consequently, it preserves biproducts. Since the codomain category
is skeletal, it follows that for any object
$(\Complex^{a_1},\ldots,\Complex^{a_n})$ in ${\SVe}^n$ we have
\begin{align*}
F(\Complex^{a_1},\ldots,\Complex^{a_n})&=
F(\Complex(1,n)^{a_1}\oplus_n\cdots\oplus_n\Complex(n,n)^{a_n})
\\ &=F(\Complex(1,n))^{a_1}\oplus_n\cdots\oplus_n F(\Complex(n,n))^{a_n}
\end{align*}
from which one easily gets that
\begin{equation} \label{accio_sobre_objectes}
F(\Complex^{a_1},\ldots,\Complex^{a_n})=(\Complex^{{\bf R}^F({\bf
a})_1},\ldots,\Complex^{{\bf R}^F({\bf a})_m})
\end{equation}
Unlike the case of linear maps, however, it
still remains to identify a minimal set of data determining the
action of $F$ on morphisms.

Let us consider the
case $m=1$. In this case, the rank matrix (or rank vector) reduces
to a point ${\bf r}=(r_1,\ldots,r_n)\in\Natural^n$, and for any
object $(\Complex^{a_1},\ldots,\Complex^{a_n})$ we have (cf.
Equation~(\ref{accio_sobre_objectes}))
\begin{equation} \label{definicio_functor_sobre_objectes}
F(\Complex^{a_1},\ldots,\Complex^{a_n})=\Complex^{{\bf a\cdot r}}
\end{equation}
For any $n\geq 1$, let us define $\Natural^n({\bf r})=\{{\bf
a}\in\Natural^n\ |\ {\bf a\cdot r}\neq 0\}$. Then, we have the
following initial description of the data defining the action
of the functor on morphisms.

%
\begin{lem} \label{lema_dades_defineixen_functor}
To give a $\Complex$-linear functor $F:{\SVe}^n\To{\SVe}$ with rank
vector ${\bf r}=(r_1,\ldots,r_n)\in\Natural^n$ it is necessary and
sufficient to give linear maps
\begin{align*}
A({\bf r};{\bf a},k,i):\Complex^{r_k}\To\Complex^{{\bf a\cdot r}}
\\
B({\bf r};{\bf a},k,i):\Complex^{{\bf a\cdot r}}\To\Complex^{r_k}
\end{align*}
for all points ${\bf a}=(a_1,\ldots,a_n)\in\Natural^n({\bf r})$,
all $k\in\{1,\ldots,n\}$ such that $r_k,a_k\neq 0$ and all
$i=1,\ldots,a_k$, satisfying the following conditions:
\begin{itemize}
\item[{\rm (A1)}]
$A({\bf r};{\bf e}_k,k',1)=B({\bf r};{\bf
e}_k,k',1)=\delta_{kk'}{\id}_{\Complex^{r_k}}$.
\item[{\rm (A2)}]
$B({\bf r};{\bf a},k,i)\circ A({\bf r};{\bf
a},k',i')=\delta_{kk'}\delta_{ii'}\ {\id}_{\Complex^{r_k}}$.
\item[{\rm (A3)}]
$\sum_{k=1}^n \sum_{i=1}^{a_k} A({\bf r};{\bf a},k,i)\circ B({\bf
r};{\bf a},k,i)={\id}_{\Complex^{{\bf a\cdot r}}}$.
\end{itemize}
(the sum over $k$ in (A3) is only over those $k$ such
that $r_k,a_k\neq 0$). Moreover, the correspondence between pairs
$(\{A({\bf r};{\bf a},k,i)\}_{{\bf a},k,i},\{B({\bf r};{\bf
a},k,i)\}_{{\bf
    a},k,i})$ as above and $\Complex$-linear functors
$F:{\SVe}^n\To{\SVe}$ with rank vector {\bf r} is a bijection.
\end{lem}
\begin{proof}
Let us consider the following presentation of ${\SVe}^n$ in terms
of $\Complex$-linear generators and relations\footnote{By
$\Complex$-linear generators of a $\Complex$-linear category $\Aa$
we obviously mean a set of non identity and non zero morphisms in
$\Aa$ such that any other morphism can be obtained from them and
from identity and zero morphisms by taking linear combinations and
compositions. As always, a relation between such a set of
generators is any pair of words in the generators defining the same
morphism in $\Aa$.}. For any $a\geq 1$ and $i=1,\ldots,a$, let
$\iota^{(a)}_i:\Complex\To\Complex^a$ (resp.
$\pi_{(a)}^i:\Complex^a\To\Complex$) be the canonical inclusion of
$\Complex$ into the $i^{th}$-component (resp. projection of
$\Complex^a$ onto the $i^{th}$-component). Furthermore, given
$a,a'\in\Natural$, denote by ${\bf
0}(a,a'):\Complex^a\To\Complex^{a'}$ the corresponding zero
morphism. Then, for any linear map $f:\Complex^a\To\Complex^{a'}$ ($a,a'\geq
1$), we have
$$
f=\sum_{i'=1}^{a'}\sum_{i=1}^a
[f]_{i'i}\ \iota^{(a')}_{i'}\circ\pi_{(a)}^i
$$
It readily follows that, given any morphism
$(f_1,\ldots,f_n):(\Complex^{a_1},\ldots,\Complex^{a_n})
\To(\Complex^{a'_1},\ldots,\Complex^{a'_n})$ in ${\SVe}^n$, it is
\begin{align*}
(f_1,\ldots,f_n)&=(f_1,{\bf 0}(a_2,a'_2),\ldots,{\bf
0}(a_n,a'_n))+\cdots+({\bf 0}(a_1,a'_1),\ldots,{\bf
0}(a_{n-1},a'_{n-1}),f_n) \\
&=\sum_{i'_1=1}^{a'_1}\sum_{i_1=1}^{a_1}[f_1]_{i'_1i_1}
(\iota^{(a'_1)}_{i'_1}\circ\pi_{(a_1)}^{i_1},{\bf
0}(a_2,a'_2),\ldots,{\bf 0}(a_n,a'_n))+\cdots \\ &\hspace{1
truecm}\cdots+
=\sum_{i'_n=1}^{a'_n}\sum_{i_n=1}^{a_n}[f_n]_{i'_ni_n}
({\bf 0}(a_1,a'_1),\ldots,{\bf 0}(a_{n-1},a'_{n-1}),
\iota^{(a'_n)}_{i'_n}\circ\pi_{(a_n)}^{i_n})
\\ &=\sum_{k=1}^n\sum_{i'_k=1}^{a'_k}\sum_{i_k=1}^{a_k}[f_k]_{i'_ki_k}
\ \iota({\bf a'},k,i'_k)\circ\pi({\bf a},k,i_k)
\end{align*}
where the maps $\iota({\bf a},k,i):({\bf
0},\ldots,\stackrel{k)}{\Complex},\ldots,{\bf
0})\To(\Complex^{a_1},\ldots,\Complex^{a_n})$ and $\pi({\bf
a},k,i):(\Complex^{a_1},\ldots,\Complex^{a_n})\To({\bf
0},\ldots,\stackrel{k)}{\Complex},\ldots,{\bf 0})$ are given by
\begin{align}
\iota({\bf a},k,i))&=({\bf
0}(0,a_1),\ldots,\stackrel{k)}{\iota^{(a_k)}_{i}},\ldots,{\bf
0}(0,a_n)) \label{generadors1} \\ \pi({\bf a},k,i)&=({\bf
0}(a_1,0),\ldots,\stackrel{k)}{\pi_{(a_k)}^{i}},\ldots,{\bf
0}(a_n,0)) \label{generadors2}
\end{align}
for any ${\bf a}=(a_1,\ldots,a_n)\in\Natural^n$ , any
$k=1,\ldots,n$ such that $a_k\neq 0$ and $i=1,\ldots,a_k$. This
shows that the maps $\{\iota({\bf a},k,i),\pi({\bf a},k,i)\}_{{\bf
a},k,i}$ are $\Complex$-linear generators for the category
${\SVe}^n$. Furthermore, we have the following obvious relations
between them:
\begin{eqnarray}
&\iota({\bf e}_k,k',1)=\pi({\bf e}_k,k',1)=\delta_{kk'}{\id}_{({\bf
0},\ldots,\stackrel{k)}{\Complex},\ldots,{\bf 0})}
\label{relacio1}
\\ &\pi({\bf a},k,i)\circ\iota({\bf
a},k',i')=\delta_{kk'}\delta_{ii'}\ {\id}_{({\bf
0},\ldots,\stackrel{k)}{\Complex},\ldots,{\bf 0})}
\label{relacio2}
\\ &\sum_{k=1}^n
\sum_{i=1}^{a_k} \iota({\bf a},k,i)\circ\pi({\bf
a},k,i)={\id}_{(\Complex^{a_1},\ldots,\Complex^{a_n})}
\label{relacio3}
\end{eqnarray}
(the sum over $k$ in (\ref{relacio3}) is only over those
values of $k$ such that $a_k\neq 0$).

Then, given a $\Complex$-linear functor $F:{\SVe}^n\To{\SVe}$ with rank vector
{\bf r}, it immediately follows from (\ref{relacio1})-(\ref{relacio3}) that
taking
\begin{align} \label{definicio_A}
A({\bf r};{\bf a},k,i)&=F(\iota({\bf a},k,i))
\\ \label{definicio_B} B({\bf r};{\bf a},k,i)&=
F(\pi({\bf a},k,i))
\end{align}
we obtain a set of data as in the statement of the Lemma.
Conversely, suppose we are given linear maps $A({\bf r};{\bf
a},k,i)$ and $B({\bf r};{\bf a},k,i)$ as above. Then, a
$\Complex$-linear functor $F:{\SVe}^n\To{\SVe}$ with rank vector
{\bf r} can be defined by
\begin{equation} \label{definicio_functor_sobre_morfismes}
F(f_1,\ldots,f_n)=\sum_{k=1}^n\sum_{i_k=1}^{a'_k}\sum_{j_k=1}^{a_k}
[f_k]_{i_kj_k}\ A({\bf r};{\bf a'},k,i_k)\circ B({\bf r};{\bf a},k,j_k)
\end{equation}
for any morphism
$(f_1,\ldots,f_n):(\Complex^{a_1},\ldots,\Complex^{a_n})\To
(\Complex^{a'_1},\ldots,\Complex^{a'_n})$ (again the sum over $k$ in
(\ref{definicio_functor_sobre_morfismes}) should be understood only
over all $k\in\{1,\ldots,n\}$
such that $a_k,a'_k\neq 0$, because otherwise $[f_k]$ is the empty
matrix). Clearly, these assignments are 
$\Complex$-linear and, by condition (A3) we have
\begin{align*}
F({\id}_{\Complex^{a_1}},\ldots,{\id}_{\Complex^{a_n}})&=
\sum_{k=1}^n\sum_{i_k=1}^{a_k}\sum_{j_k=1}^{a_k}\delta_{i_kj_k}
A({\bf r};{\bf a},k,i_k)\circ B({\bf r};{\bf a},k,j_k) \\
&=\sum_{k=1}^n\sum_{i=1}^{a_k}A({\bf r};{\bf a},k,i)\circ B({\bf
r};{\bf a},k,i)
\\ &={\id}_{\Complex^{{\bf a\cdot r}}}
\\ &={\id}_{F(\Complex^{a_1},\ldots,\Complex^{a_n})}
\end{align*}
Furthermore, given composable morphisms
$(f_1,\ldots,f_n):(\Complex^{a_1},\ldots,\Complex^{a_n})\To
(\Complex^{a'_1},\ldots,\Complex^{a'_n})$ and
$(g_1,\ldots,g_n):(\Complex^{a'_1},\ldots,\Complex^{a'_n})\To
(\Complex^{a''_1},\ldots,\Complex^{a''_n})$, condition (A2) gives that
\begin{align*}
F(g_1,\ldots,g_n)\circ F&(f_1,\ldots,f_n)=
\\
&=\left(\sum_{l=1}^n\sum_{i'_l=1}^{a''_l}
\sum_{j'_l=1}^{a'_l}[g_l]_{i'_lj'_l}\ 
A({\bf r};{\bf a''},l,i'_l)
\circ B({\bf r};{\bf a'},l,j'_l)\right)
\\ &\ \ \ \ \ \circ\left(\sum_{k=1}^n\sum_{i_k=1}^{a'_k}
\sum_{j_k=1}^{a_k}[f_k]_{i_kj_k}\ 
A({\bf r};{\bf a'},k,i_k)
\circ B({\bf r};{\bf a},k,j_k)\right)
\\ &=\sum_{k=1}^n\sum_{i'_k}^{a''_k}\sum_{j_k=1}^{a_k}
\left(\sum_{s=1}^{a'_k}[g_k]_{i'_ks}[f_k]_{sj_k}\right)
A({\bf r};{\bf a''},k,i'_k)
\circ B({\bf r};{\bf a},k,j_k) \\ &=F(g_1\circ
f_1,\ldots,g_n\circ f_n) \\ &=F((g_1,\ldots,g_n)\circ
(f_1,\ldots,f_n))
\end{align*}
so that the assignments are indeed functorial. Finally, let $F$ be
the functor defined by
(\ref{definicio_functor_sobre_objectes}) and
(\ref{definicio_functor_sobre_morfismes}) from the pair $(\{A({\bf
r};{\bf a},k,i)\}_{{\bf a},k,i},\{B({\bf r};{\bf a},k,i)\}_{{\bf
    a},k,i})$. We have
\begin{align*}
F({\bf 0}(0,a_1),\ldots,\stackrel{k)}{\iota^{(a_k)}_i},\ldots,{\bf
0}(0,a_n))&=\sum_{j=1}^{a_k}\delta_{ij} A({\bf r};{\bf a},k,j)\circ
B({\bf r};{\bf e}_k,k,1)
\\ &=A({\bf r};{\bf a},k,i)
\end{align*}
because of condition (A1), and similarly
$$
F({\bf 0}(a_1,0),\ldots,\stackrel{k)}{\pi_{(a_k)}^i},\ldots,{\bf
0}(a_n,0))=B({\bf r};{\bf a},k,i)
$$
This proves that the correspondence $(\{A({\bf r};{\bf a},k,i)\}_{{\bf
a},k,i},\{B({\bf r};{\bf a},k,i)\}_{{\bf
    a},k,i})\mapsto F$ is injective. Moreover, given $F$, it is easy to check that
if we define linear maps $A({\bf r};{\bf a},k,i)$ and $B({\bf
r};{\bf a},k,i)$ by
Equations~(\ref{definicio_A})-(\ref{definicio_B}), the functor
defined by these maps via
Equations~(\ref{definicio_functor_sobre_objectes}) and
(\ref{definicio_functor_sobre_morfismes}) is the original functor
$F$, so that the correspondence is also surjective.
\end{proof}
Although the pairs $(\{A({\bf r};{\bf a},k,i)\}_{{\bf
a},k,i},\{B({\bf r};{\bf a},k,i)\}_{{\bf a},k,i})$ already
parametrize the $\Complex$-linear functors $F:{\SVe}^n\To{\SVe}$
with rank vector {\bf r}, these data is not easy to handle. To
obtain a more manageable set of data, we proceed as follows.

Let ${\bf a}\in\Natural^n({\bf r})$ be any point such that
$a_1+\cdots+a_n\geq 2$ (namely, a point ${\bf a}\in\Natural^n({\bf
r})$ for which the $A$'s and $B$'s are not determined by condition
(A1)). Then, the corresponding two sets of linear maps $\{A({\bf
r};{\bf a},k,i)\}_{k,i}$ and $\{B({\bf r};{\bf a},k,i)\}_{k,i}$
satisfying axioms (A2) and (A3) are nothing but the data defining a
biproduct of $a_1$ copies of $\Complex^{r_1}$, $a_2$ copies of
$\Complex^{r_2}$ and so on, with object part equal to
$\Complex^{{\bf a\cdot r}}$. But in general, if
$(X,\{\iota_{X_i}\},\{\pi_{X_i}\})$ and
$(X',\{\iota'_{X_i}\},\{\pi'_{X_i}\})$ are two biproducts of the
objects $X_1,\ldots,X_n$ in any additive category $\Aa$, it is well
known that there exists a {\it unique} isomorphism $\varphi:X\To
X'$ such that
\begin{equation} \label{iso_biproducte}
\iota'_{X_i}=\varphi\circ\iota_{X_i},\quad \pi'_{X_i}=\pi_{X_i}\circ\varphi^{-1},\quad
i=1,\ldots,n
\end{equation}
and conversely, given a biproduct
$(X,\{\iota_{X_i}\},\{\pi_{X_i}\})$ and any isomorphism
$\varphi:X\To X'$, (\ref{iso_biproducte}) defines a new
biproduct. Therefore, once we have chosen a particular biproduct as
reference, given by maps $\{A^{(0)}({\bf r};{\bf a},k,i)\}_{k,i}$
and $\{B^{(0)}({\bf r};{\bf
  a},k,i)\}_{k,i}$, the two sets of maps $\{A({\bf r};{\bf a},k,i)\}_{k,i}$ and
$\{B({\bf r};{\bf
  a},k,i)\}_{k,i}$ satisfying the previous axioms are together equivalent to
a single arbitrary 
linear automorphism $\varphi({\bf a}):\Complex^{{\bf a\cdot
r}}\To\Complex^{{\bf a\cdot r}}$, related to the $A$'s and $B$'s by 
\begin{align} \label{A_en_termes_A(0)}
A({\bf r};{\bf a},k,i)&=\varphi({\bf a})\circ A^{(0)}({\bf r};{\bf
a},k,i)
\\ \label{B_en_termes_B(0)} B({\bf r};{\bf
a},k,i)&=B^{(0)}({\bf r};{\bf a},k,i)\circ \varphi({\bf a})^{-1}
\end{align}
Clearly, choosing maps $\{A^{(0)}({\bf r};{\bf a},k,i)\}_{k,i}$
and $\{B^{(0)}({\bf r};{\bf
  a},k,i)\}_{k,i}$  for all points ${\bf a}\in\Natural^n({\bf r})$
with $a_1+\cdots+a_n\geq 2$ amounts to choosing a
reference functor $H({\bf r}):{\SVe}^n\To{\SVe}$ with the given
rank vector {\bf r}. Once
such a reference functor $H({\bf r})$ has been chosen, we get a (non
canonical) bijection 
between the $\Complex$-linear functors $F:{\SVe}^n\To{\SVe}$ with
rank vector {\bf r} and the families of linear automorphisms
$$
\varphi=\{\varphi({\bf a}):\Complex^{{\bf a\cdot r}}\To\Complex^{{\bf
a\cdot r}},\ {\bf a}\in\Natural^n({\bf r})\ {\rm such\ that}\
a_1+\cdots+a_n\geq 2\}
$$
which gives the desired parametrization. By
Equations~(\ref{definicio_functor_sobre_morfismes}), (\ref{A_en_termes_A(0)})
and (\ref{B_en_termes_B(0)}), this bijection maps the
family $\varphi$ to the functor $F_{\varphi}$ acting on morphisms
by
\begin{equation} \label{definicio_functors_respecte_referencia}
F_{\varphi}(f_1,\ldots,f_n)=\varphi({\bf a'})\circ H({\bf
r})(f_1,\ldots,f_n)\circ \varphi({\bf a})^{-1}
\end{equation}
To state things more formally, let $\G({\bf r})$ be the fiber
bundle of groups over $\Natural^n$ with fibers given by
$$
\G({\bf r})_{{\bf a}}={\rm GL}({\bf a\cdot r})
$$
where ${\rm GL}(m)$ denotes the general linear group of complex
invertible $m\times m$ matrices if $m\geq 1$ and ${\rm GL}(0)=1$. Let
$\Gamma({\bf r})$ be the corresponding group 
of sections. Clearly, a section $s\in\Gamma({\bf r})$ is equivalent to a
family of automorphisms $\varphi=\{\varphi({\bf a}):\Complex^{{\bf a\cdot
r}}\To\Complex^{{\bf a\cdot r}},\ {\bf a}\in\Natural^n({\bf r})\}$, with
$\varphi({\bf a})$ such that $[\varphi({\bf a})]=s({\bf a})$ for all ${\bf
  a}\in\Natural^n({\bf r})$. Such a section will be called {\sl normalized}
when 
$$
s({\bf e}_k)=Id_{r_k}
$$
for all $k\in\{1,\ldots,n\}$ such that ${\bf e}_k\in\Natural^n({\bf
r})$. Then, we have shown the following:

\begin{thm} \label{dades_defineixen_functor}
Let ${\rm Fun}_{\Complex}^{{\bf r}}({\SVe}^n,{\SVe})$ be the set of
$\Complex$-linear functors $F:{\SVe}^n\To{\SVe}$ with rank vector
${\bf r}\in\Natural^n$, and let $\Ss({\bf r})\subset\Gamma({\bf
r})$ be the subgroup of normalized sections of the fiber bundle
$\G({\bf r})$ over $\Natural^n$. Then, there is a (non canonical)
bijection of sets
$$
{\rm Fun}_{\Complex}^{{\bf r}}({\SVe}^n,{\SVe})\cong\Ss({\bf r})
$$
\end{thm}
To be precise, for any given rank vector {\bf r}, we shall take as reference
functor $H({\bf r})$ that defined by
\begin{equation} \label{functors_referencia}
H({\bf r})(f_1,\ldots,f_n)=\bigoplus_{j=1}^n\
{\id}_{\Complex^{r_j}}\otimes f_j
\end{equation}
so that \footnote{Strictly speaking, this formula only works when
all maps $f_1,\ldots,f_n$ have nonzero domain and codomain and
$r_k\neq 0$ for all $k=1,\ldots,n$. Obviously, if $r_k=0$ and/or
$f_k$ has zero domain and codomain, the term $Id_{r_k}\otimes[f_k]$
is empty and the corresponding row and column disappear. However,
if $r_k\neq 0$ and $f_k$ has a zero domain but a nonzero codomain
(resp. a zero codomain but a nonzero domain), we should take care
of adding the appropriate number of rows (resp. columns) of zeros
in the appropriate positions, in spite of the fact that the
corresponding term $Id_{r_k}\otimes[f_k]$ is empty. In particular, this may
give rise to non diagonal matrices.}
\begin{equation} \label{matriu_functor_referencia}
[H({\bf r})(f_1,\ldots,f_n)]=\left(\begin{array}{cccc} {\rm
Id}_{r_1}\otimes [f_1] & 0 & \cdots & 0 \\ 0 & {\rm
Id}_{r_2}\otimes [f_2] & \cdots & 0 \\ \vdots & &
\ddots &
\vdots \\ 0 &
\cdots & & 0 \\ 0 & \cdots & 0 & {\rm Id}_{r_n}\otimes [f_n] \end{array}\right)
\end{equation}
The reader may easily check that this choice corresponds to taking
as maps $A^{(0)}({\bf r};{\bf a},k,i)$ and $B^{(0)}({\bf r};{\bf
  a},k,i)$ those defined by
\begin{equation} \label{A(0)}
[A^{(0)}({\bf r};{\bf a},k,i)]=[B^{(0)}({\bf r};{\bf
a},k,i)]^T=\left(
\begin{array}{c} {\bf 0}_{s,r_k}
  \\  \\ E_{i1}^{a_k,r_k} \\ \vdots \\ E_{ir_k}^{a_k,r_k} \\  \\ {\bf 0}_{t,r_k}
\end{array} \right)
\end{equation}
where ${\bf 0}_{p,q}$ denotes the $p\times q$ zero matrix,
$E_{ij}^{p,q}$ the $p\times q$ unit matrix with all entries equal
to zero except the $(i,j)$-entry, which is equal to 1, and
$s=\sum_{l=1}^{k-1}a_lr_l$, $t=\sum_{l=k+1}^na_lr_l$.

%
%

More generally, since a $\Complex$-linear functor $F:{\SVe}^n\To{\SVe}^m$ is
equivalent to $\Complex$-linear functors $F_1,\ldots,F_m:{\SVe}^n\To{\SVe}$,
for arbitrary target categories ${\SVe}^m$, $m\geq 1$, we have the following
parametrization.

\begin{cor} \label{dades_defineixen_functor_general}
Let ${\rm Fun}_{\Complex}^{{\bf R}}({\SVe}^n,{\SVe}^m)$ be the set
of $\Complex$-linear functors $F:{\SVe}^n\To{\SVe}^m$ with rank
matrix ${\bf R}\in {\rm Mat}_{m\times n}(\Natural)$. Then, there is
a (non canonical) bijection of sets
$$
{\rm Fun}_{\Complex}^{{\bf R}}({\SVe}^n,{\SVe}^m)\cong\Ss({\bf
R}_1)\times\cdots\times \Ss({\bf R}_m)
$$
\end{cor}

\begin{defn}
The $m$-tuple of normalized sections $(s_1,\ldots,s_m)\in\Ss({\bf
R}_1)\times\cdots\times \Ss({\bf R}_m)$ defining a functor
$F:{\SVe}^n\To{\SVe}^m$ with rank matrix ${\bf R}$ will be called the
{\sl gauge} of the functor (with respect to the above chosen references
$\{H({\bf r})\}_{{\bf r}}$).
\end{defn}
Explicitly, the (unique) $\Complex$-linear functor $F=(F_1,\ldots,F_m)$ with
rank 
matrix {\bf R} and gauge $(s_1,\ldots,s_m)$ is that whose
$i^{th}$-component $F_i:{\SVe}^n\To{\SVe}$, $i=1,\ldots,m$, is the
functor acting on objects and morphisms according to
Equations~(\ref{definicio_functor_sobre_objectes}) and
(\ref{definicio_functors_respecte_referencia}), respectively, with
${\bf r}={\bf R}_i$ and $\varphi_i({\bf a})$ the linear map such that
$[\varphi_i({\bf a})]=s_i({\bf a})$.

In particular, if a section $s_i$ in the gauge is {\it trivial},
i.e., such that $s_j({\bf a})=Id_{a_j}$ for all points ${\bf
a}\in\Natural^n$, the component $F_i$ of $F$ reduces to the
corresponding reference functor $H({\bf R}_i)$. This is the case,
when $n=m$, for each section in the gauge defining the identity
functor ${\id}_{{\SVe}^n}:{\SVe}^n\To{\SVe}^n$, whose rank matrix
is ${\bf R}=Id_n$. Whenever a section $s\in\Ss({\bf r})$ is trivial
in the above sense, we will write $s=1$.

Before finishing this paragraph, it is worth emphasizing that
Equation~(\ref{definicio_functors_respecte_referencia}) readily
implies the well-known fact that a $\Complex$-linear functor is
determined up to isomorphism by its rank matrix, i.e., we have the following:

\begin{cor} \label{teorema_functors_isomorfs}
Let $F,G:{\SVe}^n\To{\SVe}^m$ be two $\Complex$-linear functors
with rank matrices ${\bf R}(F)$ and ${\bf R}(G)$, respectively.
Then, they are isomorphic if and only if ${\bf R}(F)={\bf R}(G)$.
\end{cor}
Obviously, this fact already implies the existence of the correspondence
between sections of $\Gamma({\bf r})$ and $\Complex$-linear
functors with rank vector {\bf r}. Indeed, any
such functor being naturally isomorphic to
$H({\bf r})$, it is necessarily of the form
(\ref{definicio_functors_respecte_referencia}), with $\varphi({\bf
a})$ the components
of any natural isomorphism $\tau:H({\bf r})\Rightarrow F$. But it is not so
obvious from the corollary
alone that we really get a bijection by just restricting to the
normalized sections, because there will be in general various
natural isomorphisms $\tau$ between $H({\bf r})$ and $F$.

\subsection{Composition of 1-morphisms}

\label{subseccio_composicio_functors}

According to the previous subsection, after suitable
coordinatization, a $\Complex$-linear functor
$F:{\SVe}^n\To{\SVe}^m$ is completely specified by a pair $({\bf R},s)$, with
${\bf R}$ the rank matrix and $s=(s_1,\ldots,s_m)$ the gauge. In this
subsection, we find how to compute the composition of
$\Complex$-linear functors in terms of these data. This turns out
to be the more delicate point in the search of a totally
coordinatized description of {\dossv}.

Let $F:{\SVe}^n\To{\SVe}^m$ and $G:{\SVe}^m\To{\SVe}^p$ be
arbitrary $\Complex$-linear functors, defined by pairs $({\bf
R}^F,(s^F_1,\ldots,s^F_m))$ and $({\bf R}^G,(s^G_1,\ldots,s^G_p))$,
respectively. The rank matrix ${\bf R}^{G\circ F}$ of the composite
functor is easily computed. Thus, for any $j\in\{1,\ldots,n\}$, we
have (cf. (\ref{accio_sobre_C(j,n)})-(\ref{accio_sobre_objectes}))
\begin{align*}
(G\circ F)(\Complex(j,n))&=G(\Complex^{{\bf
R}^F_{1j}},\ldots,\Complex^{{\bf R}^F_{mj}})
\\ &=(\Complex^{\sum_{i=1}^m{\bf R}^G_{1i}{\bf R}^F_{ij}},\ldots,
\Complex^{\sum_{i=1}^m{\bf R}^G_{pi}{\bf R}^F_{ij}})
\end{align*}
from which it immediately follows that
$$
{\bf R}^{G\circ F}={\bf R}^G{\bf R}^F,
$$
the usual matrix product. It remains to determine how the gauge
of the composite functor $(s^{G\circ F}_1,\ldots,s^{G\circ F}_n)$ 
is computed from the gauges $(s^F_1,\ldots,s^F_m)$,
$(s^G_1,\ldots,s^G_p)$ of $F$ and $G$.

To simplify notation, we write ${\bf R}$, ${\bf R}'$ and
$\widetilde{{\bf R}}$ for the rank matrices of $F$, $G$ and $G\circ
F$, respectively, and similarly for the sections in the
corresponding gauges. By definition of the gauges, for any
$k\in\{1,\ldots,p\}$, any points ${\bf a},{\bf
a}'\in\Natural^n(\widetilde{{\bf R}}_k)$ and any morphism
$(f_1,\ldots,f_n):(\Complex^{a_1},\ldots,\Complex^{a_n})\To
(\Complex^{a'_1},\ldots,\Complex^{a'_n})$ in ${\SVe}^n$, we have
$$
\widetilde{s}_k({\bf a}')[H(\widetilde{{\bf R}}_k)(f_1,\ldots,f_n)]
\widetilde{s}_k({\bf a})^{-1}=s'_k({\bf R}({\bf a}'))
[H({\bf R}'_k)(g_1,\ldots,g_m)] s'_k({\bf R}({\bf a}))^{-1}
$$
with $g_i:\Complex^{{\bf a\cdot R}_i}\To \Complex^{{\bf a}'{\bf
\cdot R}_i}$ defined by
$$
[g_i]=s_i({\bf a}')[H({\bf R}_i) (f_1,\ldots,f_n)] s_i({\bf
a})^{-1},\quad i=1,\ldots,m
$$
Let us consider the case each $f_j$, $j=1,\ldots,n$, is an
endomorphism (hence, ${\bf a}'={\bf a}$). In this case, we know
from (\ref{matriu_functor_referencia}) that $[H(\widetilde{{\bf
R}}_k)(f_1,\ldots,f_n)]$ is the block diagonal matrix with
$\widetilde{{\bf R}}_{k1}$ copies of $[f_1]$,
$\widetilde{{\bf R}}_{k2}$ copies of $[f_2]$, etc, while
$[H({\bf R}'_k)(g_1,\ldots,g_m)]$ is the block diagonal matrix with
${\bf R}'_{k1}$ copies of $[g_1]$, ${\bf R}'_{k2}$ copies of
$[g_2]$ and so on \footnote{Notice that some of the matrices
$[f_j]$ or $[g_i]$ may be empty. However, this has no effect on the
block diagonal structure of the matrices $[H(\widetilde{{\bf
R}}_k)(f_1,\ldots,f_n)]$ and $[H({\bf R}'_k)(g_1,\ldots,g_m)]$.}.
In particular, it follows that the matrix
\begin{equation} \label{matriu_conjugada}
\left(\bigoplus_{i=1}^m Id_{{\bf R}'_{ki}}\otimes s_i({\bf a})^{-1}\right)
[H({\bf R}'_k)(g_1,\ldots,g_m)]
\left(\bigoplus_{i=1}^m Id_{{\bf R}'_{ki}}\otimes s_i({\bf a})\right)
\end{equation}
is the block diagonal matrix with ${\bf R}'_{k1}$ copies of
$[H({\bf R}_1)(f_1,\ldots,f_n)]$, ${\bf R}'_{k2}$ copies of
$[H({\bf R}_2)(f_1,\ldots,f_n)]$, etc. Now, by definition of the
functors $H({\bf R}_i)$, $i=1,\ldots,m$, this is an $m$-block
diagonal matrix, the $i^{th}$-block being itself the block
diagonal matrix made of ${\bf R}'_{ki}$ copies of
$$
{\bf diag}([f_1],\stackrel{{\bf
R}_{i1})}{\ldots},[f_1],\ldots,[f_n],\stackrel{{\bf
R}_{in})}{\ldots},[f_n])
$$
We conclude that (\ref{matriu_conjugada}) can be
obtained from $[H(\widetilde{{\bf R}}_k)(f_1,\ldots,f_n)]$ by
suitably reordering the various copies of the submatrices
$[f_1],\ldots,[f_n]$ along the diagonal. In other words, there
exist permutation matrices ${\bf P}({\bf R}'_k,{\bf R},{\bf a})$, of order
$\sum_{j=1}^n\widetilde{{\bf R}}_{kj}a_j$, such that
\begin{align} \label{condicio_matrius_permutacio}
\left(\bigoplus_{i=1}^m Id_{{\bf
R}'_{ki}}\otimes s_i({\bf a})^{-1}\right)&[H({\bf
R}'_k)(g_1,\ldots,g_m)]
\left(\bigoplus_{i=1}^m Id_{{\bf R}'_{ki}}\otimes s_i({\bf a})\right)= 
\nonumber
\\ &={\bf P}({\bf R}'_k,{\bf R},{\bf a})
[H(\widetilde{{\bf R}}_k)(f_1,\ldots,f_n)]
{\bf P}({\bf R}'_k,{\bf R},{\bf a})^{-1}
\end{align}
The annoying point is that, for given ${\bf R}'_k,{\bf R},{\bf a}$, there may
be various permutation matrices ${\bf P}({\bf R}'_k,{\bf R},{\bf a})$
for which (\ref{condicio_matrius_permutacio}) holds for all
$(f_1,\ldots,f_n):(\Complex^{a_1},\ldots,\Complex^{a_n})\To
(\Complex^{a_1},\ldots,\Complex^{a_n})$. This is because, for each
$j=1,\ldots,n$, there may be various copies of the submatrix
$[f_j]$ in $[H(\widetilde{{\bf R}}_k)(f_1,\ldots,f_n)]$. For
instance, if ${\bf a}=a{\bf e}_j$, $j=1,\ldots,n$ and $a>0$, all
maps $f_{j'}$, for $j'\neq j$, are the zero map of the
zero vector space and both (\ref{matriu_conjugada})
and $[H(\widetilde{{\bf R}}_k)(f_1,\ldots,f_n)]$ reduce to
$$
{\bf diag}([f_j],\stackrel{\widetilde{{\bf
R}}_{kj})}{\ldots},[f_j])
$$
Then, we can take as ${\bf P}({\bf R}'_k,{\bf R},a{\bf e}_j)$ any
$\widetilde{{\bf R}}_{kj}a\times\widetilde{{\bf R}}_{kj}a$ permutation matrix
of the form 
$$
{\bf P}({\bf R}'_k,{\bf R},a{\bf e}_j)={\bf P}_{\widetilde{{\bf
R}}_{kj}}\otimes Id_a
$$
with ${\bf P}_{\widetilde{{\bf R}}_{kj}}$ an arbitrary permutation
matrix of order $\widetilde{{\bf R}}_{kj}$. In particular, we can
take
\begin{equation} \label{matriu_permutacio_cas_trivial}
{\bf P}({\bf R}'_k,{\bf R},a{\bf e}_j)=Id_{\widetilde{{\bf
R}}_{kj}a}
\end{equation}
Similarly, if ${\bf R}'_k={\bf e}_i$, $i=1,\ldots,m$, both
(\ref{matriu_conjugada}) and $[H(\widetilde{{\bf
R}}_k)(f_1,\ldots,f_n)]$ reduce to
$$
{\bf diag}([f_1],\stackrel{{\bf
R}_{i1})}{\ldots},[f_1],\ldots,[f_n],\stackrel{{\bf
R}_{in})}{\ldots},[f_n])
$$
and we can take as ${\bf P}({\bf e}_i,{\bf R},{\bf a})$ any
matrix of the form
$$
{\bf P}({\bf e}_i,{\bf R},{\bf a})={\bf diag}({\bf P}_{{\bf
R}_{i1}}\otimes Id_{a_1},\ldots,{\bf P}_{{\bf R}_{in}}\otimes
Id_{a_n})
$$
with ${\bf P}_{{\bf R}_{ij}}$ an arbitrary permutation matrix of order ${\bf
R}_{ij}$, for $j=1,\ldots,n$. In particular, we can take again the
corresponding identity matrix
\begin{equation} \label{propietat_2}
{\bf P}({\bf e}_i,{\bf R},{\bf a})=Id_{\sum_{j=1}^n{\bf R}_{ij}a_j}
\end{equation}
In general, it is clear that any two possible choices ${\bf
P}({\bf R}'_k,{\bf R},{\bf a})$ and ${\bf P}({\bf R}'_k,{\bf
R},{\bf a})'$, for given ${\bf R}'_k,{\bf R},{\bf a}$, will differ by just an
arbitrary reordering, for each $j=1,\ldots,n$, of the various copies of
$[f_j]$ in $[H(\widetilde{{\bf R}}_k)(f_1,\ldots,f_n)]$. Hence, they will
necessarily be related by an equation of the form
$$
{\bf P}({\bf R}'_k,{\bf R},{\bf a})'={\bf P}({\bf R}'_k,{\bf
R},{\bf a}){\bf diag}({\bf P}_{\widetilde{{\bf R}}_{k1}}\otimes
Id_{a_1},\ldots,{\bf P}_{\widetilde{{\bf R}}_{kn}}\otimes Id_{a_n})
$$
In spite of this arbitrariness, once the
matrix ${\bf P}({\bf R}'_k,{\bf 
R},(1,\stackrel{n)}{\ldots},1))$ has been chosen, it is clear that the
matrices ${\bf P}({\bf 
R}'_k,{\bf R},{\bf a})$ for the remaining points ${\bf a}\in\Natural^n$ can
be obtained by just making them act exactly as ${\bf P}({\bf R}'_k,{\bf
R},(1,\stackrel{n)}{\ldots},1))$ except that the action is by blocks
  of dimensions $a_1,\ldots,a_n$. More explicitly, we can take as ${\bf P}({\bf
R}'_k,{\bf R},{\bf a})$ the matrix obtained after replacing each nonzero entry
  in ${\bf P}({\bf R}'_k,{\bf
R},(1,\stackrel{n)}{\ldots},1))$ by an identity or an empty matrix, according
  to the value of the corresponding $a_j$, and simultaneously deleting the
  zero entries or replacing them by suitable zero matrices
so that the resulting matrix becomes a square matrix of the
desired dimension, namely $\sum_{j=1}^n\widetilde{{\bf
R}}_{kj}a_j$.

\begin{ex} \label{exemple_matrius_permutacio} {\rm
Let ${\bf R}'_k=(1,2,1)$ and ${\bf R}=\left(\begin{array}{cc} 1 & 1 \\ 0 & 1
    \\ 2 & 0\end{array}\right)$, so that the matrix (\ref{matriu_conjugada}) is
$$
\left(\begin{array}{cccccc} [f_1]&0&0&0&0&0 \\ 0&[f_2]&0&0&0&0 \\
    0&0&[f_2]&0&0&0 
    \\ 0&0&0&[f_2]&0&0 \\ 0&0&0&0&[f_1]&0 \\ 0&0&0&0&0&[f_1]
  \end{array}\right)
$$
while $[H(\widetilde{{\bf R}}_k)(f_1,\ldots,f_n)]$ is the matrix
$$
\left(\begin{array}{cccccc} [f_1]&0&0&0&0&0 \\ 0&[f_1]&0&0&0&0 \\
    0&0&[f_1]&0&0&0 
    \\ 0&0&0&[f_2]&0&0 \\ 0&0&0&0&[f_2]&0 \\ 0&0&0&0&0&[f_2]
  \end{array}\right)
$$
Suppose we take
$$
{\bf P}((1,2,1),\left(\begin{array}{cc} 1 & 1 \\ 0 & 1 \\ 2 &
0\end{array}\right),(1,1))=\left(\begin{array}{cccccc} 1 & 0 & 0 &0
& 0 & 0
\\ 0 & 0 & 0&1 & 0 & 0 \\ 0 & 0 & 0 &0 & 1 & 0 \\ 0 & 0 & 0 &0 & 0 & 1\\
0 & 1 & 0 &0 & 0 & 0 \\ 0 & 0 & 1&0 & 0 & 0 \end{array}\right)
$$
Then, by applying the above procedure, one gets for ${\bf
  a}=(2,1)$, for example, the permutation matrix
$$
{\bf P}((1,2,1),\left(\begin{array}{cc} 1 & 1 \\ 0 & 1 \\ 2 &
0\end{array}\right),(2,1))=\left(\begin{array}{ccccccccc}
1&0&0&0&0&0&0&0&0
\\ 0&1&0&0&0&0&0&0&0 \\ 0&0&0&0&0&0&1&0&0 \\ 0&0&0&0&0&0&0&1&0 \\ 0&0&0&0&0&0&0&0&1
\\ 0&0&1&0&0&0&0&0&0 \\ 0&0&0&1&0&0&0&0&0 \\ 0&0&0&0&1&0&0&0&0 \\ 0&0&0&0&0&1&0&0&0
\end{array}\right)
$$
}
\end{ex}
To fix this ambiguity, let us give the following definition:

\begin{defn} \label{definicio_matrius_permutacio}
In the above notations, we denote by
${\bf P}({\bf R}'_k,{\bf R},{\bf 
  a})$ the matrix defined as follows:
\begin{itemize}
\item[(i)] 
${\bf P}({\bf R}'_k,{\bf
R},(1,\stackrel{n)}{\ldots},1))$ is the unique permutation matrix involving no
change in the order of the various copies of $[f_j]$, for each $j=1,\ldots,n$,
when going from $[H(\widetilde{{\bf R}}_k)(f_1,\ldots,f_n)]$ to
(\ref{matriu_conjugada}) (see example above), and
\item[(ii)]
${\bf P}({\bf R}'_k,{\bf R},{\bf a})$, for ${\bf
  a}\neq(1,\stackrel{n)}{\ldots},1)$, is that obtained
from ${\bf P}({\bf R}'_k,{\bf
R},(1,\stackrel{n)}{\ldots},1))$ by replacing each nonzero entry by an
identity or an empty matrix, according 
  to the value of the corresponding $a_j$, and simultaneously deleting the
  zero entries or replacing them by suitable zero matrices
so that the resulting matrix becomes a square matrix of dimension
$\sum_{j=1}^n\sum_{i=1}^m{\bf R}'_{ki}{\bf R}_{ij}a_j$.
\end{itemize}
\end{defn}
For what follows, it is not necessary to have an explicit expression for these
matrices. The only thing which we need is that they satisfy the following {\it
  normalization conditions}, which the reader may
easily check:

\begin{lem} \label{condicions_normalitzacio}
The matrices ${\bf P}({\bf R}'_k,{\bf R},{\bf a})$ are such that:
\begin{itemize}
\item[(i)] ${\bf P}({\bf e}_i,{\bf R},{\bf a})=Id_{\sum_{j=1}^nR_{ij}a_j}$,
  for all $i=1,\ldots,m$; 
\item[(ii)]
${\bf P}({\bf R}'_k,{\bf R},a{\bf e}_j)=
Id_{a\sum_{i=1}^m R'_{ki}R_{ij}}$, for all
$j=1,\ldots,n$;
\item[(iii)]
${\bf P}({\bf R}'_k,Id_n,{\bf
  a})=Id_{\sum_{j=1}^n R'_{kj}a_j}$, and 
\item[(iv)]
${\bf P}({\bf R}'_k,{\bf r},a)=
Id_{a\sum_{i=1}^m R'_{ki}r_i}$ when {\bf R} is an $m\times 
  1$ matrix {\bf r}, so that {\bf a} reduces to a number.
\end{itemize}
\end{lem}
Let us come back to the problem of computating the gauge of $G\circ F$. If
we define the invertible matrix
\begin{equation}
\label{calcul_seccions_2}
{\bf M}_k({\bf R}'_k,{\bf R},{\bf a})=
\widetilde{s}_k({\bf a})^{-1}s'_k({\bf R}({\bf a}))
\left(\bigoplus_{i=1}^m Id_{{\bf
R}'_{ki}}\otimes s_i({\bf a})\right){\bf P}({\bf R}'_k,{\bf R},{\bf
  a})
\end{equation}
it follows from the previous equations that
\begin{equation} \label{calcul_seccions_1}
[H(\widetilde{{\bf R}}_k)(f_1,\ldots,f_n)]={\bf M}_k({\bf
R}'_k,{\bf R},{\bf a})[H(\widetilde{{\bf
R}}_k)(f_1,\ldots,f_n)]{\bf M}_k({\bf R}'_k,{\bf R},{\bf a})^{-1}
\end{equation}
for any endomorphism
$(f_1,\ldots,f_n):(\Complex^{a_1},\ldots,\Complex^{a_n})\To
(\Complex^{a_1},\ldots,\Complex^{a_n})$ in ${\SVe}^n$. We claim that this
implies ${\bf M}_k({\bf R}'_k,{\bf R},{\bf a})$ is of the form
\begin{equation} \label{calcul_seccions_4}
{\bf M}_k({\bf R}'_k,{\bf R},{\bf a})=\bigoplus_{j=1}^n\
\Lambda_j({\bf
  R}'_k,{\bf R},{\bf a})\otimes Id_{a_j}
\end{equation}
for some invertible $\widetilde{{\bf R}}_{kj}\times \widetilde{{\bf
R}}_{kj}$ matrices $\Lambda_j({\bf
  R}'_k,{\bf R},{\bf a})$, $j=1,\ldots,n$. Indeed, let ${\bf M}_k({\bf
  R}'_k,{\bf R},{\bf a})$ 
be the matrix
$$
\left(\begin{array}{ccccccc} B^{11}_{11} & \cdots & B^{11}_{1\widetilde{{\bf R}}_{k1}} & \cdots &
B^{1n}_{11} & \cdots & B^{1n}_{1\widetilde{{\bf R}}_{kn}} \\
\vdots & & \vdots & & \vdots & & \vdots \\
B^{11}_{\widetilde{{\bf R}}_{k1}1} & \cdots &
B^{11}_{\widetilde{{\bf R}}_{k1}\widetilde{{\bf R}}_{k1}} &
\cdots & B^{1n}_{\widetilde{{\bf R}}_{k1}1} & \cdots & B^{1n}_{\widetilde{{\bf R}}_{k1}
\widetilde{{\bf R}}_{kn}} \\ \vdots & & \vdots & & \vdots & & \vdots \\
B^{n1}_{11} & \cdots & B^{n1}_{1\widetilde{{\bf R}}_{k1}} &
\cdots & B^{nn}_{11} & \cdots & B^{nn}_{1\widetilde{{\bf
R}}_{kn}} \\
\vdots & & \vdots & & \vdots & & \vdots \\
B^{n1}_{\widetilde{{\bf R}}_{kn}1} & \cdots &
B^{n1}_{\widetilde{{\bf R}}_{kn}\widetilde{{\bf R}}_{k1}} &
\cdots & B^{nn}_{\widetilde{{\bf R}}_{kn}1} & \cdots & B^{nn}_{\widetilde{{\bf R}}_{kn}
\widetilde{{\bf R}}_{kn}} \end{array} \right)
$$
with $B^{jj'}_{rs}$ an $a_j\times a_{j'}$ complex matrix,
$j,j'=1,\ldots,n$. Then, (\ref{calcul_seccions_1}) says that this
matrix commutes with all block diagonal matrices of the form
$$
{\bf diag}(A_1,\stackrel{\widetilde{{\bf
R}}_{k1})}{\ldots},A_1,\ldots,A_n,\stackrel{\widetilde{{\bf
R}}_{kn})}{\ldots},A_n)
$$
for arbitrary $A_1,\ldots,A_n$, with $A_j$ an $a_j\times a_j$
complex matrix. It readily follows that $B^{jj'}_{rs}=0$ if $j\neq
j'$ and that $B^{jj}_{rs}=\lambda^{(j)}_{rs}({\bf a})\ Id_{a_j}$
with $\lambda^{(j)}_{rs}({\bf a})\in\Complex$. Then, the matrix
$\Lambda_j({\bf R}'_k,{\bf R},{\bf a})$ referred to in
(\ref{calcul_seccions_4}) is that whose entries are the numbers
$\lambda^{(j)}_{rs}({\bf a})$. Since ${\bf M}_k({\bf R}'_k,{\bf
R},{\bf a})$ is invertible, all these matrices $\Lambda_j({\bf
R}'_k,{\bf R},{\bf a})$ will also be invertible.

In fact, the coefficients $\lambda^{(j)}_{rs}({\bf a})$ are
independent of {\bf a}. This is easily seen using that
(\ref{calcul_seccions_1}) actually holds for any morphism
$(f_1,\ldots,f_n):(\Complex^{a_1}\ldots,\Complex^{a_n})
\To(\Complex^{a'_1}\ldots,\Complex^{a'_n})$ in ${\SVe}^n$,
not necessarily an endomorphism, a fact which we leave as an
exercise to the reader to prove. Thus, when $(f_1,\ldots,f_n)$ is
not an endomorphism, each matrix $A_j$ above is an arbitrary
rectangular $a_j\times a'_j$ matrix and condition
(\ref{calcul_seccions_1}) establishes that
$$
\lambda^{(j)}_{rs}({\bf a}')A_j=\lambda^{(j)}_{rs}({\bf a})A_j
$$
for all $j,r,s$. Hence, we can write
\begin{equation} \label{calcul_seccions_3}
{\bf M}_k({\bf R}'_k,{\bf R},{\bf a})=
\bigoplus_{j=1}^n\
\Lambda_j({\bf
  R}'_k,{\bf R})\otimes Id_{a_j}
\end{equation}
Now, the matrices $\Lambda_j({\bf R}'_k,{\bf R})$ can be easily computed using
that all sections in any gauge are normalized. Thus, letting ${\bf
a}={\bf e}_j$ in (\ref{calcul_seccions_2}) and
(\ref{calcul_seccions_3}), we get that
$$
\Lambda_j({\bf
  R}'_k,{\bf R})={\bf M}_k({\bf R}'_k,{\bf R},{\bf e}_j)=
s'_k({\bf R}({\bf e}_j)){\bf P}({\bf R}'_k,{\bf R},{\bf
  e}_j)=s'_k({\bf R}({\bf e}_j))
$$
because of the second normalization condition
(cf. Lemma~\ref{condicions_normalitzacio}). Consequently 
$$
{\bf M}_k({\bf R}'_k,{\bf R},{\bf a})=\bigoplus_{j=1}^n\ s'_k({\bf
R}({\bf e}_j))\otimes Id_{a_j}
$$
and, coming back to (\ref{calcul_seccions_2}), we obtain for
$\widetilde{s}_k({\bf a})$ the expression:
$$
\widetilde{s}_k({\bf a})=s'_k({\bf R}({\bf a}))
\left(\bigoplus_{i=1}^m Id_{{\bf
R}'_{ki}}\otimes s_i({\bf a})\right){\bf P}({\bf R}'_k,{\bf R},{\bf
  a})
\left(\bigoplus_{j=1}^n\ s'_k({\bf R}({\bf e}_j))^{-1}\otimes
Id_{a_j}\right)
$$
Therefore, we have proved the following rule, expressing the
composition of $\Complex$-linear functors in terms of the pairs $({\bf R},s)$: 

\begin{thm}
Let $F:{\SVe}^n\To{\SVe}^m$ and $G:{\SVe}^m\To{\SVe}^p$ be
arbitrary $\Complex$-linear functors, defined by rank matrices
${\bf R}^F$ and ${\bf R}^G$ and gauges $(s^F_1,\ldots,s^F_m)$ and
$(s^G_1,\ldots,s^G_p)$, respectively. Then, the rank matrix and
gauge of the composite functor $G\circ F$ are given by
\begin{equation} \label{matriu_functor_composicio}
{\bf R}^{G\circ F}={\bf R}^G{\bf R}^F
\end{equation}
and
\begin{align} \label{seccions_functor_composicio}
s_k^{G\circ F}({\bf a})=s_k^G({\bf R}^F({\bf a}))&
\left(\bigoplus_{i=1}^m Id_{{\bf R}^G_{ki}}\otimes
s_i^F({\bf a})\right){\bf P}({\bf R}^G_k,{\bf R}^F,{\bf
  a}) \nonumber \\ &\hspace{1.3
truecm}\left(\bigoplus_{j=1}^n\ s^G_k({\bf R}^F({\bf
e}_j))^{-1}\otimes Id_{a_j}\right)
\end{align}
for all $k=1,\ldots,p$.
\end{thm}
Notice that the sections defined by
(\ref{seccions_functor_composicio}) are indeed normalized because of the
second normalization condition on the permutation matrices ${\bf P}({\bf R}^G_k,{\bf R}^F,{\bf a})$. Furthermore, when one of the functors $F$
or $G$ is an identity (hence, the corresponding sections are
trivial), Equation~(\ref{seccions_functor_composicio}) clearly
gives the sections of the other functor, so that $(Id_n,1)$ really acts as a
unit for this 
composition law between pairs $({\bf R},s)$. Less obvious is the strict
associativity of this composition, but in our
approach this follows from the 
associativity of functor composition.

\subsection{2-morphisms}

As done previously with 1-morphisms, in this subsection we
determine a set that parametrizes the 2-morphism in
{\dossv} between two given 1-morphisms.

Let $F,F':{\SVe}^n\To{\SVe}^m$ be arbitrary $\Complex$-linear
functors, with rank matrices {\bf R} and ${\bf R}'$ respectively,
and let $\tau:F\Rightarrow F'$ be a natural transformation. A
priori, giving $\tau$ requires giving all its components
$$
\tau_{(\Complex^{a_1},\ldots,\Complex^{a_n})}:F(\Complex^{a_1},\ldots,\Complex^{a_n})\To
F'(\Complex^{a_1},\ldots,\Complex^{a_n})
$$
for all objects $(\Complex^{a_1},\ldots,\Complex^{a_n})$ in
${\SVe}^n$. However, by using $k$-order biproduct functors
$\oplus_n^{(k)}:{\SVe}^n\times\stackrel{k)}{\cdots}\times{\SVe}^n\To{\SVe}^n$
on ${\SVe}^n$ for each $k\geq 2$, it turns out that a few of these
components allows one to compute the remaining ones. Indeed, an
arbitrary object 
$(\Complex^{a_1},\ldots,\Complex^{a_n})$ in ${\SVe}^n$ is the
object part of a $(\sum_{j=1}^na_j)$-order biproduct, namely
\begin{equation} \label{objectes_igual_biproductes}
(\Complex^{a_1},\ldots,\Complex^{a_n})=\oplus_n^{(\sum_{j=1}^n
a_n)}(\Complex(1,n),\stackrel{a_1)}{\ldots},\Complex(1,n),\ldots,\Complex(n,n),
\stackrel{a_n)}{\ldots},\Complex(n,n)),
\end{equation}
and we have the following well-known general fact:

\begin{lem} \label{lemaa}
Let $\Aa,\Bb$ be categories, with $\Aa$ additive and $\Bb$
preadditive. Suppose given two additive functors $F,G:\Aa\To\Bb$
and let $\tau:F\Rightarrow G$ be a natural transformation. Then,
for any objects $A_1,\ldots,A_k$ of $\Aa$, we have
\begin{equation} \label{tau_suma_directa}
\tau_{A_1\oplus\cdots\oplus A_k}=\sum_{i=1}^k G(\iota_{A_i})\circ\tau_{A_i}
\circ F(\pi_{A_i})
\end{equation}
where $A_1\oplus\cdots\oplus A_k$, $\iota_{A_i}:A_i\To
A_1\oplus\cdots\oplus A_k$ and $\pi_{A_i}:A_1\oplus\cdots\oplus
A_k\To A_i$, $i=1,\ldots,k$, denote the object and the morphisms
defining {\sl any} biproduct of $A_1,\ldots,A_k$.
\end{lem}
\begin{proof}
By naturality of $\tau$, we have
$G(\pi_{A_i})\circ\tau_{A_1\oplus\cdots\oplus A_k}=\tau_{A_i}\circ
F(\pi_{A_i})$ for all $i=1,\ldots,k$. Hence, taking the composite
on the left with $G(\iota_{A_i})$ we get
$$
G(\iota_{A_i}\circ\pi_{A_i})\circ\tau_{A_1\oplus\cdots\oplus A_k}=
G(\iota_{A_i})\circ\tau_{A_i}\circ F(\pi_{A_i})
$$
Equation~(\ref{tau_suma_directa}) then follows by suming up over
all $i=1,\ldots,k$ and using that $\sum_{i=1}^k
\iota_{A_i}\circ\pi_{A_i}={\id}_{A_1\oplus\cdots\oplus A_k}$.
\end{proof}

Notice from the previous Lemma that we can choose any biproduct of
$A_1,\ldots,A_k$ to compute the $(A_1\oplus\cdots\oplus A_k)$-component
of $\tau$. Then, we shall fix $k$-order biproduct functors $\oplus_n^{(k)}$
such that, for any ${\bf a}=(a_1,\ldots,a_n)\in\Natural^n$, the
biproduct (\ref{objectes_igual_biproductes}) is that for
which the defining morphisms
$\iota^{j,l_j}_{a_1,\ldots,a_n}:\Complex(j,n)\To(\Complex^{a_1},\ldots,\Complex^{a_n})$
and
$\pi^{j,l_j}_{a_1,\ldots,a_n}:(\Complex^{a_1},\ldots,\Complex^{a_n})\To\Complex(j,n)$,
for all $l_j=1,\ldots,a_j$ and all $j=1,\ldots,n$, are given by
\begin{align} \label{iota_biproductes_basics}
\iota^{j,l_j}_{a_1,\ldots,a_n}
&=\iota({\bf a}, j,l_j)
\\ \label{pi_biproductes_basics} \pi^{j,l_j}_{a_1,\ldots,a_n}
&=\pi({\bf a}, j,l_j)
\end{align}
with $\iota({\bf a}, j,l_j)$ and $\pi({\bf a}, j,l_j)$ the
linear generators of ${\SVe}^n$ introduced in
Lemma~\ref{dades_defineixen_functor}. 

It follows from (\ref{tau_suma_directa}) that $\tau$ is uniquely
determined by the components
$\{\tau_{\Complex(j,n)},\ j=1,\ldots,n\}$. Furthermore, the reader
may easily check that these components can be chosen arbitrarily, because
naturality is always fulfilled when the remaining components are
defined by (\ref{tau_suma_directa}).

According to (\ref{accio_sobre_C(j,n)}),
$\tau_{\Complex(j,n)}:F(\Complex(j,n))\To F'(\Complex(j,n))$ is of the form
$$
\tau_{\Complex(j,n)}=(\Phi(\tau)_{1j}\ldots,\Phi(\tau)_{mj})
$$
for some linear maps $\Phi(\tau)_{ij}:\Complex^{R_{ij}}\To\Complex^{R'_{ij}}$
for all $i=1,\ldots,m$. Fixing canonical bases in each 
vector space (the last step in the coordinatization process), this
is further equivalent to a family of complex matrices ${\bf T}(\tau)_{ij}$
(some of them possibly empty). Therefore,
associated to $\tau$, we have an
$m\times n$ matrix ${\bf T}(\tau)$ whose $(i,j)$-entry is in turn a
complex matrix of type $R'_{ij}\times R_{ij}$ if $R'_{ij},R_{ij}\neq 0$ and
empty otherwise. Such a 
matrix ${\bf T}(\tau)$ will be called the {\sl matrix of} $\tau$
and we have proved the following:

\begin{thm} \label{identificacio_2-morfismes}
Let $F,F':{\SVe}^n\To{\SVe}^m$ be $\Complex$-linear functors, with
rank matrices {\bf R} and {\bf R}', and let ${\rm Mat}_{m\times
n}({\bf R},{\bf R}')$ be the set of $m\times n$ matrices {\bf T}
whose $(i,j)$-entry is an $R'_{ij}\times R_{ij}$ complex
matrix if $R'_{ij},R_{ij}\neq 0$ and empty otherwise. Then, there is a
(non canonical) bijection of sets
$$
{\rm Nat}(F,F')\cong{\rm Mat}_{m\times n}({\bf R},{\bf R}')
$$
mapping a natural transformation $\tau$ to its matrix ${\bf
T}(\tau)$.
\end{thm}
Obviously, the identity natural transformation of any functor is mapped by
this bijection to the trivial matrix {\bf T} with all non empty entries equal
to identities.

Observe that the parametrizing set ${\rm Mat}_{m\times
  n}({\bf 
R},{\bf R}')$ in the previous bijection is the same for all pairs of functors
$F,F'$ with the 
same rank matrices {\bf R} and ${\bf R}'$ (in agreement with the
fact that a $\Complex$-linear functor is determined up to
isomorphism by its rank matrix), so that an
element ${\bf T}\in{\rm Mat}_{m\times n}({\bf R},{\bf R}')$
represents various natural transformations, insofar as the domain and codomain
functors 
$F$ and $F'$ are not specified among all functors in the respective
isomorphism classes.

Let us finally remark that Theorem~\ref{identificacio_2-morfismes} could also
be 
deduced using an enriched version of
Yoneda's lemma. Thus, it is easily checked that a $\Complex$-linear
functor $F:{\bf SVect}^n\To{\bf SVect}$ of rank ${\bf
  r}=(r_1,\ldots,r_n)$ can be represented by the object
$(\Complex^{r_1},\ldots,\Complex^{r_n})$, so that, if $G:{\bf SVect}^n\To{\bf SVect}$ is of rank ${\bf s}=(s_1,\ldots,s_n)$, we have bijections
$$
{\rm Nat}(F,G)\cong{\rm Nat}({\rm Hom}_{{\bf
    SVect}^n}((\Complex^{r_1},\ldots,\Complex^{r_n}),-),G)\cong
G(\Complex^{r_1},\ldots,\Complex^{r_n})\cong\Complex^{r_1s_1+\ldots+r_ns_n}
$$
and this can indeed be
identified with the set ${\rm Mat}_{1\times n}({\bf r},{\bf
  s})$. However, these identifications are also non canonical. Our
previous approach has the advantage
that it makes easier to find how compositions between 2-morphisms
should be computed in terms of the corresponding matrices. This is done in the
next paragraph.

\subsection{Vertical and horizontal composition of 2-morphisms}

Let $F,F',F'':{\SVe}^n\To{\SVe}^m$ be $\Complex$-linear functors,
respectively defined by pairs $({\bf R}^F,s^F)$,
$({\bf R}^{F'},s^{F'})$, $({\bf
R}^{F''},s^{F''})$, and let $\tau:F\Rightarrow
F'$ and $\tau':F'\Rightarrow F''$ be natural transformations,
described by matrices ${\bf T}(\tau)$ and ${\bf T}(\tau')$. It is
very easy to see how to compute the matrix ${\bf
T}(\tau'\cdot\tau)$ of the vertical composite of $\tau$ and
$\tau'$.

\begin{thm}
In the above notations, if $R^{F''}_{ij},R^{F'}_{ij},R^F_{ij}\neq 0$ it is
\begin{equation} \label{composicio_vertical_b}
{\bf T}(\tau'\cdot\tau)_{ij}={\bf T}(\tau')_{ij}\ {\bf
T}(\tau)_{ij},\quad i=1,\ldots,m,\ j=1,\ldots,m
\end{equation}
the usual matrix product. Otherwise, it is ${\bf
  T}(\tau'\cdot\tau)_{ij}=\emptyset$ or ${\bf
  T}(\tau'\cdot\tau)_{ij}={\bf 0}$ (a zero matrix), depending on which
  coefficients $R^{F''}_{ij},R^{F'}_{ij},R^F_{ij}$ are zero. 
\end{thm}
\begin{proof}
The formula readly follows from the definitions of ${\bf T}(\tau)$
and ${\bf T}(\tau')$ and the fact that, for any
$j\in\{1,\ldots,n\}$, it is
$(\tau'\cdot\tau)_{\Complex(j,n)}=\tau'_{\Complex(j,n)}\circ\tau_{\Complex(j,n)}$.
\end{proof}
It is worth noting that, according to (\ref{composicio_vertical_b}),
${\bf T}(\tau'\cdot\tau)$ is independent of the sections defining
the various functors involved. It depends only on the matrices
${\bf T}(\tau)$ and ${\bf T}(\tau')$.

This is no longer true as regards horizontal composition. Indeed, let
us consider additional $\Complex$-linear functors 
$G,G':{\Ve}^m\To{\Ve}^p$, described by pairs $({\bf
R}^G,s^G)$ and $({\bf
R}^{G'},s^{G'})$, and a natural transformation
$\sigma:G\Rightarrow G'$, described by a matrix ${\bf T}(\sigma)$.
We know that
$$
(\sigma\circ\tau)_{\Complex(j,n)}=\sigma_{F'(\Complex(j,n))}\circ
G(\tau_{\Complex(j,n)}),\quad j=1,\ldots,n
$$
Now, for all $j=1,\ldots,n$, it is
\begin{align*}
F'(\Complex(j,n))&=(\Complex^{R^{F'}_{1j}},\ldots,\Complex^{R^{F'}_{mj}}) \\
&=\oplus_m^{(\sum_{i=1}^m
  R^{F'}_{ij})}(\Complex(1,m),\stackrel{R^{F'}_{1j})}{\ldots},
\Complex(1,m),\ldots,\Complex(m,m),\stackrel{R^{F'}_{mj})}{\ldots},
\Complex(m,m)) 
\end{align*}
Hence, by Equation~(\ref{tau_suma_directa}) and using the above higher order
biproduct functors (cf. 
Equations~(\ref{iota_biproductes_basics})-(\ref{pi_biproductes_basics})),
the $k$-component of $\sigma_{F'(\Complex(j,n))}$, for any
$k=1,\ldots,p$, is given by (cf. (\ref{definicio_A})-(\ref{definicio_B}))
\begin{align*}
\left[(\sigma_{F'(\Complex(j,n))})_k\right]&=\sum_{i=1}^m\sum_{l_i=1}^{{\bf
R}^{F'}_{ij}} \left[G'_k(\iota({\bf R}^{F'}({\bf
e}_j),i,l_i))\right](\sigma_{\Complex(i,m)})_k\left[G_k(\pi({\bf
R}^{F'}({\bf e}_j),i,l_i))\right]
\\ &=\sum_{i=1}^m\sum_{l_i=1}^{{\bf
R}^{F'}_{ij}} \left[A^{G'}_k({\bf R}^{G'}_k;{\bf R}^{F'}({\bf
e}_j),i,l_i)\right]{\bf T}(\sigma)_{ki}\left[B^G_k({\bf
R}^{G}_k;{\bf R}^{F'}({\bf e}_j),i,l_i)\right]
\\ &=s^{G'}_k({\bf R}^{F'}({\bf
e}_j))\Xi(F',{\bf T}(\sigma),k,j) s^G_k({\bf R}^{F'}({\bf
e}_j))^{-1}
\end{align*}
where (cf.
(\ref{A_en_termes_A(0)})-(\ref{B_en_termes_B(0)}))
$$
\Xi(F',{\bf T}(\sigma),k,j)=
\sum_{i=1}^m\sum_{l_i=1}^{{\bf R}^{F'}_{ij}}
\left[A^{(0)}({\bf R}^{G'}_k;{\bf R}^{F'}({\bf e}_j),i,l_i)\right]{\bf
T}(\sigma)_{ki}\left[B^{(0)}({\bf R}^{G}_k;{\bf R}^{F'}({\bf
e}_j),i,l_i)\right]
$$
(see the conventions adopted in the introduction as regards
notation). Moreover, the $k^{th}$ component
$G_k(\tau_{\Complex(j,n)})$ of the term 
$G(\tau_{\Complex(j,n)})$, for any $k=1,\ldots,p$, is given by (cf.
(\ref{definicio_functors_respecte_referencia}))
\begin{align*}
\left[G_k(\tau_{\Complex(j,n)})\right]&=
\left[G_k(\Phi(\tau)_{1j},\ldots,\Phi(\tau)_{mj})\right] 
\\ &=s^G_k({\bf R}^{F'}({\bf e}_j))\left[H({\bf 
R}^G_k)(\Phi(\tau)_{1j},\ldots,\Phi(\tau)_{mj})\right]
s^G_k({\bf R}^F({\bf e}_j))^{-1}
\end{align*}
Therefore
\begin{align*}
{\bf T}(\sigma\circ\tau)_{kj}&=\left[(\sigma_{F'(\Complex(j,n))})_k\right]
\left[G_k(\tau_{\Complex(j,n)})\right] \\ &=s^{G'}_k({\bf R}^{F'}({\bf
e}_j))\ \Xi(F',{\bf T}(\sigma),k,j) \\ &\hspace{1 truecm}
\left[H({\bf R}^G_k)(\Phi(\tau)_{1j},\ldots,\Phi(\tau)_{mj})\right]
s^G_k({\bf R}^F({\bf e}_j))^{-1}
\end{align*}
for all $j=1,\ldots,n$, $k=1,\ldots,p$.

\begin{lem}
For all $j=1,\ldots,n$, $k=1,\ldots,p$, it is
$$
\Xi(F',{\bf T}(\sigma),k,j)\left[
H({\bf
  R}^G_k)(\Phi(\tau)_{1j},\ldots,\Phi(\tau)_{mj})\right]=\bigoplus_{i=1}^m\
{\bf T}(\sigma)_{ki}\otimes{\bf 
T}(\tau)_{jk}
$$
(recall that ${\bf T}(\sigma)_{ki}$ denotes the matrix in canonical basis of
the linear map $\Phi(\sigma)_{ki}$, and similarly for ${\bf T}(\tau)_{jk}$). 
\end{lem}
\begin{proof}
Since $B^{(0)}({\bf R}^{G}_k;{\bf R}^{F'}({\bf e}_j),i,l_i)=H({\bf
R}^G_k)(\pi({\bf R}^{F'}({\bf e}_j),i,l_i))$, we have
\begin{align*}
B^{(0)}({\bf R}^{G}_k;{\bf R}^{F'}({\bf e}_j),i,l_i)&\circ H({\bf
R}^G_k)(\Phi(\tau)_{1j},\ldots,\Phi(\tau)_{mj})=
\\ &=H({\bf R}^G_k)
({\bf 0}(R^F_{1j},0),\ldots,\stackrel{i)}{\pi^{l_i}_{(R^{F'}_{ij})}\circ
\Phi(\tau)_{ij}},\ldots,{\bf 0}(R^F_{mj},0))
\end{align*}
for any $i,l_i$. It suffices then to see that, for all
$i=1,\ldots,m$, the matrix
\begin{align*}
\sum_{l_i=1}^{R^{F'}_{ij}}&
\left[A^{(0)}({\bf R}^{G'}_k;{\bf
R}^{F'}({\bf e}_j),i,l_i)\right]{\bf
T}(\sigma)_{ki} \\ &\hspace{1 truecm}\left[ H({\bf R}^G_k)({\bf 0}
(R^F_{1j},0),\ldots,\stackrel{i)}{\pi^{l_i}_{
(R^{F'}_{ij})}\circ\Phi(\tau)_{ij}},\ldots,{\bf 0}(R^F_{mj},0))\right]
\end{align*}
is zero everywhere except for an $(R^{G'}_{ki}R^{F'}_{ij})
\times(R^{G}_{ki}R^{F}_{ij})$ non zero block, equal to
${\bf T}(\sigma)_{ki}\otimes{\bf T}(\tau)_{ij}$ and starting at
row $1+s_i=1+\sum_{q=1}^{i-1} R^{G'}_{kq}R^{F'}_{qj}$
and column $1+t_i=1+\sum_{q=1}^{i-1} R^{G}_{kq}R^{F}_{qj}$. Now, it easily
follows from (\ref{A(0)}) that 
$[A^{(0)}({\bf R}^{G'}_k;{\bf R}^{F'}({\bf
  e}_j),i,l_i)\circ\Phi(\sigma)_{ki}]$ looks like 
$$
\left(\begin{array}{c} {\bf 0}_{s_i\times R^G_{ki}} \\ \\
{\bf M}({\bf T}(\sigma)_{ki};l_i,1)
\\ \\ \hspace{0.5 truecm}\vdots\ \ {\bf R}^{G'}_{ki}) \\ \\
{\bf M}({\bf T}(\sigma)_{ki};l_i,{\bf R}^{G'}_{ki}) \\
\\ {\bf 0}_{\tilde{s}_i\times {\bf R}^G_{ki}}
\end{array}\right)
$$
where ${\bf M}({\bf T}(\sigma)_{ki};l_i,q))$, $q=1,\ldots,
R^{G'}_{ki}$, is the $R^{F'}_{ij}\times R^G_{ki}$
matrix everywhere zero except for the $l_i^{th}$ row, which is
equal to the $q^{th}$ row of ${\bf T}(\sigma)_{ki}$, while
$\tilde{s}_i=\sum_{q=i+1}^m R^{G'}_{kq}R^{F'}_{qj}$.
Similarly, by definition of the reference functors $H$, we have
that $[H({\bf R}^G_k)({\bf 0}(R^F_{1j},0),\ldots,\stackrel{i)}{\pi^{l_i}_{(
R^{F'}_{ij})}\circ\Phi(\tau)_{ij}},\ldots,{\bf 0}(R^F_{mj},0))]$ looks like
$$
\left(\begin{array}{ccc} {\bf 0}_{R^G_{ki}\times t_i} & {\bf N}({\bf T}(\tau);l_i)
& {\bf 0}_{R^G_{ki}\times \tilde{t}_i}
\end{array}\right)
$$
where ${\bf N}({\bf T}(\tau)_{ij};l)$ denotes the block diagonal
matrix
$$
{\bf N}({\bf T}(\tau)_{ij};l)={\rm diag}(\mbox{$l^{th}$-row\ of}\
{\bf T}(\tau)_{ij},\ \stackrel{R^G_{ki})}{\ldots}\
,\mbox{$l^{th}$-row\ of}\ {\bf T}(\tau)_{ij})
$$
while $\tilde{t}_i=\sum_{q=i+1}^m R^{G}_{kq}R^{F}_{qj}$. The above statement
readily follows now by taking the 
product of both matrices and suming up over all $l_i=1,\ldots,R^{F'}_{ij}$.
\end{proof}
Hence, we have proved the following:

\begin{thm}
For any natural transformations $\tau:F\Rightarrow F'$ and
$\sigma:G\Rightarrow G'$, the matrix ${\bf T}(\sigma\circ\tau)$ of the
horizontal composite is given by
\begin{equation} \label{composicio_horitzontal_b}
{\bf T}(\sigma\circ\tau)_{kj}=s^{G'}_k({\bf R}^{F'}({\bf
e}_j))\ \left(\bigoplus_{i=1}^m\ {\bf T}(\sigma)_{ki}\otimes{\bf
T}(\tau)_{jk}\right)\ s^G_k({\bf R}^F({\bf e}_j))^{-1}
\end{equation}
for all $j=1,\ldots,n$ and $k=1,\ldots,p$.
\end{thm}
Except for the gauge terms $s^{G'}_k({\bf R}^{F'}({\bf
e}_j))$ and $s^{G}_k({\bf R}^{F}({\bf
e}_j))$, this is the formula that defines the horizontal
composition between 2-morphisms in Kapranov and Voevodsky's totally
coordinatized version ${\dosv}_{cc}$. We leave as an exercise to the
reader to directly check that, unlike Kapranov and Voevodsky's
formula, (\ref{composicio_horitzontal_b}) indeed
defines a strictly associative composition. Notice also that,
according to this formula, the identity natural transformation of any
identity functor act as a unit with respect to horizontal composition,
as required in any strict 2-category.

\section{The strict 2-category ${\dossv}_{cc}$}

The previous arguments allow us to formulate the definition of
our strict totally 
coordinatized version of {\dosv} as follows.

\begin{defn}
${\bf 2SVect}_{cc}$ is the
strict 2-category with objects, 1-morphisms and 2-morphisms defined by:

\begin{itemize}
\item
objects: the natural numbers $n\geq 0$.
\item
1-morphisms: given $n,m\geq 1$, a 1-morphism $n\To m$ is a pair
$({\bf R},s)$, where
\begin{itemize}
\item[(i)]
${\bf R}$ (the {\sl rank matrix}) is an $m\times n$ matrix $(R_{ij})$
with entries in $\Natural$, and
\item[(ii)]
$s$ (the {\sl gauge}) is a collection $\{(s_1({\bf a}),\ldots,s_m({\bf
  a}))\}_{{\bf a}\in\Natural^n}$ with $s_i({\bf a})$ an ${\bf
  R}({\bf a})_i\times {\bf R}({\bf a})_i$ 
invertible complex matrix if ${\bf R}({\bf a})_i\neq 0$ and
$s_i({\bf a})=1$ otherwise, and satisfying the following {\sl normalization
condition}: for all $i=1,\ldots,m$ and $j=1,\ldots,n$, it is
\begin{equation} \label{condicio_normalitzacio_gauge}
s_i({\bf e}_j)=Id_{R_{ij}}
\end{equation}
whenever $R_{ij}={\bf R}({\bf e}_j)_i\neq 0$ (otherwise, it is equal to 1).
\end{itemize}
When $n=0$ and/or $m=0$, there
is only one 1-morphism $n\To m$, which will be denoted by ${\bf
0}_{n,m}$.
\item
2-morphism: given objects $n,m\geq 1$ and 1-morphisms $({\bf
  R},s),({\bf R}',s'):n\To
m$, with ${\bf R}=(R_{ij})_{i=1,\ldots,m,j=1,\ldots,n}$ and ${\bf
  R}'=(R'_{ij})_{i=1,\ldots,m,j=1,\ldots,n}$, a 2-morphism $({\bf
  R},s)\Rightarrow({\bf R}',s')$ is an
$m\times n$ matrix {\sf T} whose $(i,j)$ entry ${\sf T}_{ij}$ is an
$R'_{ij}\times R_{ij}$ complex matrix if both
$R_{ij},R'_{ij}\neq 0$ and empty otherwise. If $n=0$ and/or $m=0$,
both 1-morphisms are necessarily equal and there is a unique
2-morphism between them denoted $1_{{\bf 0}_{n,m}}$ (obviously, it
is the identity 2-morphism of ${\bf 0}_{n,m}$).
\end{itemize}
The various compositions are as follows:
\begin{itemize}
\item
composition of 1-morphisms: if $n,m,p\geq 1$ and $({\bf R},s):n\To m$ and
$(\widetilde{{\bf
  R}},\widetilde{s}):m\To p$, the composite
$(\widetilde{{\bf R}},\widetilde{s})\circ({\bf R},s):n\To p$ is the
pair
\begin{equation} \label{composicio_1-morfismes}
(\widetilde{{\bf R}},\widetilde{s})\circ({\bf
R},s)=(\widetilde{{\bf R}}{\bf R},\widetilde{s}\ast s)
\end{equation}
where $\widetilde{{\bf R}}{\bf R}$ denotes the usual matrix product
and $\widetilde{s}\ast s$ is defined by
\begin{align} \label{seccions_functor_composicio_b}
(\widetilde{s}\ast s)_k({\bf a})=\widetilde{s}_k({\bf R}({\bf a}))&
\left(\bigoplus_{i=1}^m Id_{\widetilde{R}_{ki}}\otimes
s_i({\bf a})\right){\bf P}(\widetilde{{\bf R}}_k,{\bf R},{\bf
  a}) \nonumber \\ &\hspace{1.3
truecm}\left(\bigoplus_{j=1}^n\ \widetilde{s}_k({\bf R}({\bf
e}_j))^{-1}\otimes Id_{a_j}\right)
\end{align}
for all ${\bf a}\in\Natural^n$ and $k=1,\ldots,p$ ($Id_0$
means the empty matrix). Here,
${\bf P}(\widetilde{{\bf R}}_k,{\bf R},{\bf a})$ denote the
permutation matrices of order
$\sum_{j=1}^n\sum_{i=1}^m\widetilde{R}_{ki}R_{ij}a_j$ introduced in
Definition~\ref{definicio_matrius_permutacio} and  
satisfying the normalization conditions in
Lemma~\ref{condicions_normalitzacio}. In case one of the numbers
$n,m,p$ is zero, the composite is the corresponding zero map.
\item
vertical composition of 2-morphisms: given 2-morphisms ${\sf
T}:({\bf
  R},s)\Rightarrow({\bf R}',s')$ and
${\sf T}':({\bf R}',s')\Rightarrow({\bf R}'',s'')$, with $({\bf
R},s),({\bf R}',s'),({\bf R}'',s''):n\To m$ and $n,m\geq 1$, the
vertical composite ${\sf T}'\cdot{\sf T}:({\bf R},s)\Rightarrow({\bf
  R}'',s'')$ is the matrix obtained by 
componentwise multiplication,
i.e.
\begin{equation} \label{composicio_vertical}
({\sf T}'\cdot{\sf T})_{ij}={\sf T}'_{ij}{\sf T}_{ij}, \quad
    i=1,\ldots,m,\ j=1,\ldots,n
\end{equation}
where we agree that the product of a matrix by the empty matrix
is the empty matrix or the appropriate zero matrix (for example,
if $R'_{ij}=0$ but $R_{ij},R''_{ij}\neq 0$, we have ${\sf T}_{ij}={\sf
  T}'_{ij}=\emptyset$, while $({\sf T}'\cdot{\sf T})_{ij}$
is the $R''_{ij}\times R_{ij}$ zero matrix).

\item
horizontal composition of 2-morphisms: given 2-morphisms ${\sf
T}:({\bf
  R},s)\Rightarrow({\bf R}',s'):n\To m$ and $\widetilde{{\sf T}}:(\widetilde{{\bf R}},\widetilde{s})\Rightarrow(\widetilde{{\bf R}}',\widetilde{s}'):m\To
p$, with $n,m,p\geq 1$, the horizontal composite $\widetilde{{\sf
T}}\circ{\sf T}$ is the $p\times n$ matrix
with entries defined by
\begin{equation} \label{composicio_horitzontal}
(\widetilde{{\sf T}}\circ{\sf T})_{kj}=\widetilde{s}'_k({\bf
R}'({\bf e}_j))\left(\bigoplus_{i=1}^m\ \widetilde{{\sf
T}}_{ki}\otimes{\sf T}_{ij}\right)\widetilde{s}_k({\bf R}({\bf
e}_j))^{-1}
\end{equation}
for all $j=1,\ldots,n$ and $k=1,\ldots,p$, where we agree again
that the tensor product of any matrix by the empty matrix is the
empty matrix.
\end{itemize}
\end{defn} 
The reader may easily check from the
previous formulas that the pair 
$(Id_n,1)$, where $1$ denotes the trivial gauge, indeed corresponds to the
identity 1-morphism ${\id}_n$, for any $n\geq 1$, while it follows from 
(\ref{composicio_vertical}) that the identity 2-morphisms correspond
to matrices all of 
whose nonempty entries are identities, i.e., if
$({\bf R},s):n\To m$, ${\bf
1}_{({\bf R},s)}$ is given by $({\bf 1}_{({\bf
  R},s)})_{ij}=Id_{R_{ij}}$ if $R_{ij}\neq 0$ and empty
otherwise.

It also follows immediately from the previous formulas that
a 1-morphism $({\bf R},s):n\To m$ in ${\bf 2SVect}_{cc}$ is invertible
if and only if 
$n=m$ and its rank matrix {\bf R} is a permutation
matrix, and that two 1-morphisms $({\bf R},s)$ and $({\bf R}',s')$ are
2-isomorphic if and only if 
${\bf R}={\bf R}'$ (cf. Corollary~\ref{teorema_functors_isomorfs}), a
2-morphism ${\sf T}:({\bf 
R},s)\Rightarrow({\bf R},s')$ being a 2-isomorphism if and only if all
nonempty entries in {\sf T} are non singular complex matrices.

At this point, the obvious next step is to try to materialize the
monoidal structure 
on ${\dossv}_{cc}$ inherited from the monoidal structure on
{\dosv}. This is a particularly important point to be able to carry
out the study of the representation 
theory of categorical groups on 2-vector spaces, because it is this
monoidal structure 
that induces the monoidal structure on the 2-category of
representations of the categorical group. However, making explicit
this structure is deferred to another paper.

\vspace{0.7 truecm}
\noindent{{\bf Acknowledgements.}} It is usual at this point to
ackowledge some mathematical colleagues or some institution who has
given 
some sort of financial support. Let me thank, however, the person
who has really made possible this work: my wife Merc\`e Serra.

\bibliographystyle{plain}
\bibliography{stc_KV_2cat}

\end{document}